\documentclass[12pt]{amsart2000}
\usepackage{amssymb}

\begin{document}
\numberwithin{equation}{section}
\def\Label#1{\label{#1}}

\def\1#1{\ov{#1}}
\def\2#1{\widetilde{#1}}
\def\3#1{\mathcal{#1}}
\def\4#1{\widehat{#1}}
\def\5#1{\mathbb{#1}}

\def\s{s}
\def\eps{\varepsilon}
\def\k{\kappa}
\def\ov{\overline}
\def\span{\text{\rm span}}
\def\tr{\text{\rm tr}}
\def\GL{{\sf GL}}
\def\U{{\sf U}}
\def\SU{{\sf SU}}
\def\xo {{x_0}}
\def\Rk{\text{\rm Rk\,}}
\def\sg{\sigma}
\def \emxy{E_{(M,M')}(X,Y)}
\def \semxy{\scrE_{(M,M')}(X,Y)}
\def \jkxy {J^k(X,Y)}
\def \gkxy {G^k(X,Y)}
\def \exy {E(X,Y)}
\def \sexy{\scrE(X,Y)}
\def \hn {holomorphically nondegenerate}
\def\hyp{hypersurface}
\def\prt#1{{\partial \over\partial #1}}
\def\det{{\text{\rm det}}}
\def\wob{{w\over B(z)}}
\def\co{\chi_1}
\def\po{p_0}
\def\fb {\bar f}
\def\gb {\bar g}
\def\Fb {\ov F}
\def\Gb {\ov G}
\def\Hb {\ov H}
\def\zb {\bar z}
\def\wb {\bar w}
\def \qb {\bar Q}
\def \t {\tau}
\def\z{\chi}
\def\w{\tau}
\def\Z{\zeta}

\def \T {\theta}
\def \Th {\Theta}
\def \L {\Lambda}
\def\b{\beta}
\def\a{\alpha}
\def\o{\omega}
\def\l{\lambda}

\def \im{\text{\rm Im }}
\def \re{\text{\rm Re }}
\def \Char{\text{\rm Char }}
\def \supp{\text{\rm supp }}
\def \codim{\text{\rm codim }}
\def \Ht{\text{\rm ht }}
\def \Dt{\text{\rm dt }}
\def \hO{\widehat{\mathcal O}}
\def \cl{\text{\rm cl }}
\def \bR{\mathbb R}
\def \bC{\mathbb C}
\def \bP{\mathbb P}
\def \C{\mathbb C}
\def \R{\mathbb R}
\def \bL{\mathbb L}
\def \bZ{\mathbb Z}
\def \bN{\mathbb N}
\def \scrF{\mathcal F}
\def \scrK{\mathcal K}
\def \scrM{\mathcal M}
\def \cR{\mathcal R}
\def \scrJ{\mathcal J}
\def \scrA{\mathcal A}
\def \scrO{\mathcal O}
\def \scrV{\mathcal V}
\def \scrL{\mathcal L}
\def \scrE{\mathcal E}
\def \hol{\text{\rm hol}}
\def \aut{\text{\rm aut}}
\def \Aut{\text{\rm Aut}}
\def \J{\text{\rm Jac}}
\def\jet#1#2{J^{#1}_{#2}}
\def\gp#1{G^{#1}}
\def\gpo{\gp {2k_0}_0}
\def\emmp {\scrF(M,p;M',p')}
\def\rk{\text{\rm rk}}
\def\Orb{\text{\rm Orb\,}}
\def\Exp{\text{\rm Exp\,}}
\def\Span{\text{\rm span\,}}
\def\d{\partial}
\def\D{\3J}
\def\pr{{\rm pr}}
\def\dbl{[\![}
\def\dbr{]\!]}
\def\nl{|\!|}
\def\nr{|\!|}

\def \D{\text{\rm Der}\,}
\def \Rk{\text{\rm Rk}\,}
\def \ima{\text{\rm im}\,}

\title[Rigidity of CR-immersions]{Rigidity of CR-immersions into spheres}
\author[P. Ebenfelt, X. Huang, D. Zaitsev]
{Peter Ebenfelt, Xiaojun Huang, and Dmitri Zaitsev}
\footnotetext {\rm The first
author is a Royal Swedish Academy of Sciences Research Fellow
supported in part by a grant from the Knut and Alice Wallenberg
Foundation and by DMS-0100110}
\footnotetext {\rm The second author is supported in
part by DMS-9970439 and a grant from the Rutgers Research Council}
\footnotetext {\rm The third author is supported in
part by a grant from the Italian Consiglio Nazionale delle Ricerche}
\address{P. Ebenfelt: Department of Mathematics, University
of California at San Diego, La Jolla, CA 92093, USA}
\email{pebenfel@math.ucsd.edu }
\address{X. Huang: Department of Mathematics, Rutgers University,
New Brunswick, NJ 08502, USA} \email{huangx@math.rutgers.edu}
\address{D. Zaitsev: Dipartimento di Matematica,
Universit\`a di Padova, via G. Belzoni 7, 35131 Padova, ITALY}
\email{zaitsev@math.unipd.it}
\newtheorem{Thm}{Theorem}[section]
\newtheorem{Cor}[Thm]{Corollary}
\newtheorem{Pro}[Thm]{Proposition}
\newtheorem{Lem}[Thm]{Lemma}
\theoremstyle{definition}
\newtheorem{Def}[Thm]{Definition}
\newtheorem{Rem}[Thm]{Remark}
\newtheorem{Ex}[Thm]{Example}

\begin{abstract}
We consider local CR-immersions of a strictly pseudoconvex real
hypersurface $M\subset\bC^{n+1}$, near a point $p\in M$, into the
unit sphere $\mathbb S\subset\bC^{n+d+1}$ with $d>0$. Our main
result is that if there is such an immersion $f\colon (M,p)\to
\mathbb S$ and $d < n/2$, then $f$ is {\em rigid} in the sense that
any other immersion of $(M,p)$ into $\mathbb S$ is of the form
$\phi\circ f$, where $\phi$ is a biholomorphic automorphism of the
unit ball $\mathbb B\subset\bC^{n+d+1}$. As an application of this
result, we show that an isolated singularity of an irreducible
analytic variety of codimension $d$ in $\bC^{n+d+1}$ is uniquely
determined up to affine linear transformations by the local CR
geometry at a point of its Milnor link.
\end{abstract}

\maketitle

\section{Introduction}

Let $X$ be a complex-analytic $(n+1)$-dimensional (not necessarily
closed) irreducible variety in $\C^{n+d+1}$ and $\5S_\eps$ a real
hypersphere of radius $\eps>0$ such that the open ball $\mathbb
B_\eps$ (whose boundary is $\5S_\eps$) as well as $\mathbb S_\eps$
intersect $X$; for convenience, we shall assume that the ball
$\mathbb B_\eps$ is centered at $0$. The intersection
$K_\eps:=X\cap \5S_\eps$ is then a real-analytic variety (whose
set of regular points is a real hypersurface) in $X$. Let us
denote by $M_\eps\subset K_\eps$ the relatively open subset of
points where $X$ is nonsingular and meets the sphere
transversally.  For instance, if $X$ has an isolated singularity
at $0\in X$, then, for sufficiently small generic $\eps>0$,
$M_\eps$ is all of $K_\eps$; in this case, $M_\eps$ is sometimes
called the {\em link of the singularity $(X,0)$}. In what follows,
we shall refer to $M_\eps$ as the $\mathbb S_\eps$-link (or simply
the link) of the variety $X$, even though we do not necessarily
assume that $\eps$ is small or that $X$ has an isolated
singularity at $0$ (or even that the central point $0$ is on $X$).

If $X'$ is another variety of dimension $n+1$ in $\bC^{n+d+1}$ and
$M'_\eps$ is the link of $X'$, and if there is a biholomorphic
automorphism of the ball $\5B_\eps$ which sends the variety $X$ to
$X'$, then clearly the links $M_\eps$ and $M'_\eps$ are
CR-equivalent submanifolds of $\5S_\eps$; we should point out here
that the automorphisms of the ball $\mathbb B_\eps$ are linear
fractional transformations of $\bC^{n+d+1}$, which yield, by
restriction, all the CR automorphisms of the sphere $\mathbb
S_\eps$. One may ask the converse:
\medskip

\noindent (Q) {\em Does the local CR-structure of the manifold
$M_\eps$ determine the variety $X$ uniquely $($up to biholomorphic
automorphisms of $\5B_\eps)$?}
\medskip

That is, if $X'$ is another irreducible analytic variety of the
same dimension such that the $\mathbb S_\eps$-link of $X'$ is
locally CR-equivalent to that of $X$, must then the varieties $X$
and $X'$ be equivalent by a biholomorphic automorphism of the ball
$\mathbb B_\eps$? The answer is in general ``no'' as can be seen
by the following well known example.

\begin{Ex}\Label{ex-whitney}
Consider the so-called Whitney embedding $W\colon\bC^{n+1}\to
\bC^{2n+1}$,
\begin{equation}\Label{eq-whitney}
W(z_1,\ldots, z_{n+1}) :=
(z_1,\ldots,z_{n},z_1z_{n+1},z_2z_{n+1},\ldots,
z_{n}z_{n+1},z_{n+1}^2),
\end{equation}
which is two-to-one and hence the image is an irreducible
analytic variety $X\subset\bC^{2n+1}$ (of codimension $d=n$).
The variety $X$ is singular along the the complex line
$z_1=\cdots=z_n=w_1=\ldots=w_{n}=0$, where
$(z_1,\ldots,z_{n},w_1,\ldots,w_{n+1})$ denote the coordinates in
$\bC^{2n+1}$. It is well known and easy to check that $W$ maps the
unit sphere in $\bC^{n+1}$ into the unit sphere in $\bC^{2n+1}$.
Consequently, the link $M_1$ of $X$ is {\em locally spherical}
(i.e. locally CR-equivalent to the sphere), and hence locally
equivalent to $M'_1$, where $M'_1$ is the intersection of the unit
sphere with $(n+1)$-plane
$X'=\{(z_1,\ldots,z_{n},w_1,\ldots,w_{n+1})\colon
z_1=\ldots=z_{n}=0\}$. However, the singular variety $X$ is
clearly not biholomorphically equivalent (in any neighborhood of
$0$) to the nonsingular variety $X'$.
\end{Ex}
\medskip

On the other hand, it follows from the work of Webster
\cite{We79}, Faran~\cite{Faran86}, Cima-Suffridge~\cite{CS83} that
if the link of $X$ is locally spherical and if the codimension $d$
is {\it strictly less} than $n$ (to exclude the situation in
Example~\ref{ex-whitney} above) then $X$ must be an $(n+1)$-plane.
This answers
 in the affirmative the
question (Q) above in this special case.
See also the work of Cima-Suffridge ~\cite {CS90},
 Forstneri\v c \cite{Forst}, and the second author \cite{Hu1} in some other
 related directions.

In this paper, we shall consider the general form of the question
(Q). We shall prove that if the codimension of the variety $X$ in
$\bC^{n+d+1}$ is less than half its dimension minus one, then the
answer to (Q) is affirmative. More generally, we consider the
situation where $M_\eps=M$ is a smooth (meaning here always
$\3C^\infty$) abstract CR-manifold of hypersurface type (a {\em
CR-hypersurface}) and prove the following:

\begin{Thm}\Label{thm-main0}
Let $M$ be a connected smooth CR-hypersurface of dimension $2n+1$.
If $d<n/2$, then any smooth CR-immersion $f$ of $M$ into the unit
sphere $\5S$ in $\C^{n+d+1}$ is rigid. That is, any other smooth
CR-immersion $\tilde f\colon M\to \5S$ is related to $f$ by
$\tilde f=\phi\circ f$, where $\phi$ is a CR-automorphism of
$\5S$.
\end{Thm}

In the case $d=1$, the conclusion in Theorem \ref{thm-main0}
follows from the work of Webster \cite{We79}.

Theorem~\ref{thm-main0} can be directly applied to study the
question (Q) above. Let us first consider the case where $X$ has
an isolated singularity at $0\in X$ and no other singularity of
$X$ is inside the ball $\mathbb B_\eps$. If $X'$ is another such
variety and if their links $M_\eps$, $M'_\eps$ are locally CR
equivalent at $p\in M_\eps$ and $p'\in M'_\eps$, then in view of
Theorem~\ref{thm-main0} (the unit sphere in that theorem can of
course be replaced by a sphere of any radius) there is a
biholomorphic automorphism $\phi$ of $\mathbb B_\eps$ (extending
as a birational transformation of the whole space $\bC^{n+d+1}$)
which sends an open piece of $M_\eps$ to an open piece of
$M'_\eps$. If we, in addition, assume that $X\cap \mathbb B_\eps$
and $X'\cap \mathbb B_\eps$ are connected, then we conclude that
$\phi$ sends $X\cap \mathbb B_\eps$ to $X'\cap \mathbb B_\eps$.
Since the only singularities of $X$ and $X'$ inside $\mathbb
B_\eps$ are at 0, $\phi$ must also send 0 to 0. By also using the
fact that the only biholomorphic automorphisms of the ball
$\mathbb B_\eps$ which fix the origin are {\it unitary linear
transformation} of $\bC^{n+d+1}$, we obtain the following:

\begin{Cor}\Label{thm-isolated}
Let $X, X'\subset \bC^{n+d+1}$ be irreducible complex analytic
varieties of codimension $d$ through $0$. Let $M_\eps$, $M'_\eps$
denote their manifolds of regular intersection with a sphere
$\5S_\eps$ of  radius $\eps>0$ centered at $0\in \bC^{n+d+1}$.
Assume that both $X$ and $X'$ have isolated singularities at $0$,
no other singular points in $\mathbb B_\eps$, and $X\cap\mathbb
B_\eps$, $X'\cap\mathbb B_\eps$ are connected. If $d<n/2$ and
$M_\eps$ and $M'_\eps$ are locally CR-equivalent at some points
$q\in M_\eps$ and $q'\in M'_\eps$, then there exists a unitary
linear transformation which maps $X\cap \5B_\eps$ to $X'\cap
\5B'_\eps$.
\end{Cor}

We can also consider the general case where $X$ and $X'$ are not
assumed to have isolated singularities (and $0$ is not necessarily
a point on $X\cap X'$). A similar argument to the one above yields
the following result:

\begin{Cor}\Label{thm-variety1}
Let $X, X'\subset \bC^{n+d+1}$ be irreducible complex analytic
varieties of codimension $d$ and let $M_\eps$, $M'_\eps$ denote
their manifolds of regular intersection with a sphere $\5S_\eps$
of radius $\eps>0$ centered at $0\in \bC^{n+d+1}$. If $d<n/2$ and
$M_\eps$ and $M'_\eps$ are locally CR-equivalent at some points
$q\in M_\eps$ and $q'\in M'_\eps$, then there exists a
biholomorphic automorphism of the ball $\mathbb B_\eps$ whose
birational extension to $\bC^{n+d+1}$ sends $(X,q)$ to $(X',q')$
$($and $X$ into a subvariety containing $X'$ as an open subset$)$.
\end{Cor}

In a different direction, we also obtain the following result as a
direct consequence of Theorem \ref{thm-main0}; the simple argument
needed to derive it from Theorem \ref{thm-main0} is left to the
reader.

\begin{Cor}\Label{thm-global} Let $D$ be a connected relatively compact
open subset of a
complex manifold and assume that $D$ has a smooth strictly
pseudoconvex boundary $\partial D$. Assume that there is a proper
embedding of $D$ into the unit ball $\mathbb B_1\subset
\bC^{n+d+1}$, with $d<n/2$, which extends smoothly to the boundary
$\partial D$. If $g$ is a local CR-diffeomorphism sending a
relatively open piece of $\partial D$ into $\partial D$, then $g$
extends as a biholomorphic automorphism of $D$.
\end{Cor}

Theorem \ref{thm-main0} above will follow from the more general
Theorem~\ref{thm-main3} below, where a higher codimension $d$ is
allowed provided that $M_\eps$ is suitably degenerate. As a
corollary, we recover the result about $(n+1)$-planes mentioned
above.

The proof of Theorem \ref{thm-main3}
will be completed in \S\ref{sec-deg}. It decomposes
naturally into two parts. The first part consists of showing that
if the mapping is degenerate (in a certain sense to be defined in
the next section), then the image $f(M)$ is in fact contained in
the intersection of the sphere with a lower dimensional
complex plane (see Theorem \ref{thm-main2}).
Using this fact, we reduce the proof
to the case where $f$ is finitely nondegenerate (see \S\ref{sec-main}).
The reader is referred to \S\ref{sec-main} for the statements
of mentioned results and
a more detailed outline of the proof of Theorem \ref{thm-main3}.

\section{Further results}\Label{sec-main}

Let $M$ be CR-hypersurface of dimension
$2n+1$ and denote by $\3V=T^{0,1}M\subset\C TM$ its CR-bundle; the
reader is referred e.g. to \cite{BER} for basic notions and facts
about CR-manifolds (see also \S\ref{sec-prelim}).
Recall that a mapping $f=(f_1,\ldots,f_k)\colon M\to \bC^k$ is
called CR if $f_*(\3V_p)\subset T_{f(p)}^{0,1}\bC^k$ for every $p\in M$.
This is equivalent to saying that $L f_j=0$ for every $j=1,\ldots k$
and every $(0,1)$ vector field $L$ (i.e.\ a section in $\3V$).

An important step in proving Theorem~\ref{thm-main0} will be to
classify possible CR-immersions according to their degeneracy. Let
$\hat M\subset \bC^{\hat n+1}$ be a real hypersurface (and hence
also a CR-manifold) and $f\colon M\to \bC^{\hat n+1}$ be a CR
mapping sending $M$ into $\hat M$. We shall refer to $d=\hat n-n$
as the codimension of the mapping $f$; thus, $\hat M$ is a real
hypersurface in $\bC^{n+d+1}$ . Let $p$ be a point in $M$ and
$\hat\rho$ a local defining function for $\hat M$ near $\hat
p:=f(p)\in \hat M$. Following Lamel \cite{Lamel}, define an
increasing sequence of subspaces $E_k(p)$ of the space $T'_{\hat
p}\C^{n+d+1}$ of $(1,0)$ covectors as follows. Let $L_{\bar
1},\ldots, L_{\bar n}$ be a basis of $(0,1)$ vector fields on $M$
near $p$ and define
\begin{equation}\Label{eq-ksdeg1}
E_k(p):=\span_\C\{L^{\bar J}(\hat \rho_{Z'}\circ f)(p) : J\in
(\mathbb Z_+)^n, 0\le|J|\le k\} \subset  T'_{\hat p}\C^{n+d+1},
\end{equation}
where $\hat \rho_{Z'}=\d\hat\rho$ is represented by vectors in
$\C^{n+d+1}$ in some local coordinate system $Z'$ near $\hat p$;
we use here standard multi-index notation $ L^{\bar J}:=L_{\bar
1}^{\bar J_1}\ldots L_{\bar n}^{\bar J_n}$ and
$|J|:=J_1+\ldots+J_n$. One can show (cf.\ \cite{Lamel}) that
$E_k(p)$ is independent of the choice of local defining function
$\hat \rho$ and coordinates $Z'$, as well as of the choice of
basis of the CR vector fields $L_{\bar 1},\ldots, L_{\bar n}$. We
shall say, again following Lamel \cite{Lamel,L2} that $f\colon
(M,p)\to \hat M$ is {\em $(k_0,s_0)$-degenerate} at $p$ if $s_0$
is the minimum of the decreasing sequence of codimensions of
$E_k(p)$ in $\C^{n+d+1}$, i.e.\
\begin{equation}\Label{deg1}
s_0=s(p):=\min_k \codim E_k(p)
\end{equation}
and $k=k_0$ is the smallest integer for which this minimum is
attained. If $E_k(p)=T'_{\hat p}\C^{n+d+1}$ for some $k$, i.e.\ if
$f$ is $(k_0,0)$-degenerate at $p$ for some $k_0$, it is said to
be {\em $k_0$-nondegenerate} (at $p$) or {\em finitely
nondegenerate} without specifying $k_0$.

Let us call the {\em degeneracy} of $f$ the minimum of $s(p)$ as
defined in \eqref{deg1} for $p\in M$. We have the following result
which, in view of Theorem \ref{thm-secform2} and Proposition
\ref{pro-basic} below, implies Theorem \ref{thm-main0}.

\begin{Thm}\Label{thm-main3}
Let $f\colon M\to \5S$ be a smooth CR-immersion of a smooth
connected CR-hypersurface $M$ of dimension $2n+1$ into the unit
sphere $\5S$ in $\bC^{n+d+1}$ and denote by $s$ be the degeneracy
of $f$. If $d-s<n/2$, then $f$ is rigid among smooth CR-immersions
having the same degeneracy. That is, any other such CR-immersion
$\tilde f\colon M\to \5S$ is related to $f$ by $\tilde f=\phi\circ
f$, where $\phi$ is a CR-automorphism of $\5S$.
\end{Thm}

It is not difficult to see from the definition that, if $f(M)$ is
contained in a complex plane of codimension $s$ in $\C^{n+d+1}$,
then the degeneracy of $f$ is at least $s$. A important ingredient
in the proof of Theorem \ref{thm-main3} is the following (partial)
converse, which also seems to be of independent interest.

\begin{Thm}\Label{thm-main2}
Let $f\colon M\to \5S$ be a smooth CR-immersion of a smooth
connected CR-hypersurface $M$ of dimension $2n+1$ into the unit
sphere $\5S$ in $\bC^{n+d+1}$. Let $s$ be the degeneracy of $f$.
If $d-s<n$, then $f(M)$ is contained in the intersection of $\5S$
with a complex plane $P\subset \bC^{n+d+1}$ of codimension $s$.
Moreover, if $f\colon M\to \5S$ is $(k_0,s)$-nondegenerate at a
point $p\in M$, then it is $k_0$-nondegenerate there as a mapping
$f\colon M\to \5S\cap P$.
\end{Thm}

Another important ingredient in the proof of
Theorem~\ref{thm-main0} will be the CR analogue of the classical
Gauss equation which, as the reader may recall, relates the
Riemannian curvature tensors of a manifold and its submanifold
with the second fundamental form composed with the Riemannian
metric (see e.g.\ \cite{KN2}). The {\em $($extrinsic$)$ second
fundamental form} for a CR-mapping $f\colon M\to \hat M$ between
real hypersurfaces $M\subset\C^{n+1}$ and $\hat
M\subset\C^{n+d+1}$
 can be defined (up to a scalar factor) by
\begin{equation}\Label{eq-defsec}
\Pi(X_{p},Y_{p}):=\overline{\pi\big(XY{(\hat \rho_{\bar Z'}\circ
f)}(p)\big)} \in \1{T'_{\hat p}\hat M/E_1(p)},
\end{equation}
where $\pi\colon T'_{\hat p}\hat M\to T'_{\hat p}\hat M/E_1(p)$
is the projection and $X$, $Y$ are any $(1,0)$
vector fields on $M$ extending given vectors $X_{p},Y_{p}\in T^{1,0}_pM$.
In the case where $\hat M$ (and hence also $M$) is strictly pseudoconvex,
the Levi form of $\hat M$ (at $\hat p$) with respect to $\hat\rho$ defines
an isomorphism
$\1{T'_{\hat p}\hat M/E_1(p)}\cong T^{1,0}_{\hat p}\hat M/f_* T^{1,0}_pM$
and hence the second fundamental form can be viewed
as an $\C$-bilinear symmetric form
\begin{equation}\Label{SFF}
\Pi_p\colon T^{1,0}_pM\times T^{1,0}_pM\to T^{1,0}_{\hat p}\hat M/f_* T^{1,0}_pM
\end{equation}
that does not depend anymore on the choice of $\hat\rho$. In
\S\ref{sec-subs} (see equation (\ref{eq-secid})), we shall show
that $\Pi$ is indeed the second fundamental form (in the classical
sense) of $f$ with respect to Webster's pseudohermitian connection
induced by the CR-structure (and a choice of contact form). We
shall say that the second fundamental form of $f$ is {\em
nondegenerate at $p$} if its values span the target space. This is
easily seen to be equivalent to $E_2(p)=T'_{\hat p} \C^{n+d+1}$
or, in the notation above, to $f$ being $2$-nondegenerate at $p$.

The CR analogue of the Weyl curvature tensor in Riemannian
geometry is given by the tangential pseudoconformal curvature
tensor
$$S\colon T^{1,0}_pM\times T^{1,0}_pM\times T^{1,0}_pM\times T^{1,0}_pM\to
\C T_pM/(T^{1,0}_pM\oplus T^{0,1}_pM)$$ defined by Chern and Moser
\cite{CM} for every Levi-nondegenerate CR-hypersurface $M$ (in
fact, it is defined there as a tensor on a principal bundle over
$M$ which can be pulled back to $M$ as will be explained in
\S\ref{sec-prelim}, see also \cite{We79}). The role of the
Riemannian metric itself is played by the Levi form which can be
invariantly seen as a Hermitian bilinear map
$$L\colon T^{1,0}_pM\times T^{1,0}_pM\to \C T_pM/(T^{1,0}_pM\oplus T^{0,1}_pM).$$
In the following we identify the quotient $T^{1,0}_{\hat p}\hat
M/f_* T^{1,0}_pM$ with the orthogonal complement of $f_*
T^{1,0}_pM$ in $T^{1,0}_{\hat p}\hat M$ with respect to the Levi
form of $\hat M$. We shall also use the notation $f_*$ for the
mapping $\C T_pM/(T^{1,0}_pM\oplus T^{0,1}_pM)\to \C \hat T_{\hat
p}M/(T^{1,0}_{\hat p}\hat M\oplus T^{0,1}_{\hat p}\hat M)$ induced
by $f_*\colon \bC T_p M\to \bC T_{\hat p}\hat M$. The CR analogue
of the Gauss identity mentioned above can be now stated as
follows.

\begin{Thm}\Label{gauss-main}
Let $f\colon M\to \hat M$ be a smooth CR-immersion of
a smooth connected CR-hypersurface $M$
into a strongly pseudoconvex CR-hypersurface $\hat M$.
Denote by $\Pi$ the second fundamental form of $f$,
by $L$ and $S$ the Levi form and the tangential pseudoconformal curvature tensor
for $M$ and by $\hat L$ and $\hat S$ the corresponding tensors for $\hat M$.
Then there exists a Hermitian form $H\colon T^{1,0}_pM\times T^{1,0}_pM\to \C$
such that the identity
\begin{multline}\Label{psc-gauss1}
\hat S(f_*V,f_*V,f_*V,f_*V) - f_* S(V,V,V,V)  = \hat
L\big(\Pi(V,V), \Pi(V,V)\big)  \\ + f_* L(V,V)\, H(V,V)
\end{multline}
holds for every  $V \in T^{1,0}_pM$.
\end{Thm}

The reader is referred to \S\ref{sec-gauss} for more details on
this Gauss equation and to Proposition~\ref{psconf-gauss} from
which it immediately follows. Here we only briefly mention that
the term involving $H$ on the right hand side (which will be
called conformally equivalent to $0$, cf. \S\ref{sec-prelim}) can
be expressed from \eqref{psc-gauss1} in terms of $\hat S$, $\hat
L$, $L$ and $\Pi$ and hence \eqref{psc-gauss1} can be given a more
explicit form (see \eqref{gauss-explicit}).

We conclude this section by giving an outline of the
paper and proofs of Theorems \ref{thm-main2} and \ref{thm-main3}.
In \S\ref{sec-prelim} we recall the construction of the pseudohermitian connection
defined by Webster~\cite{We78} on a given strongly pseudoconvex CR-hypersurface
with a fixed contact form and show the relation to the pseudoconformal connection
defined by Chern and Moser~\cite{CM} in terms of the forms \eqref{eq-CRpar}.
In \S\ref{sec-subs} we show the existence
of coframes suitably adapted to a pair $(M,\hat M)$ of strongly pseudoconvex CR-hypersurfaces,
where $M$ is a submanifold of $\hat M$. It is further shown the relation
of the pseudohermitian connection with respect to such a coframe with the second
fundamental form as defined above. \S\ref{sec-gauss} is devoted to
the proof of pseudohermitian and pseudoconformal analogues of the Gauss equation
(Propositions~\ref{psherm-gauss} and \ref{psconf-gauss}
from which Theorem~\ref{gauss-main} follows).
As one of the main consequences (see Corollary~\ref{thm-secform}(ii))
we obtain the following, which seems interesting in its
own right:

{\em Under the assumptions of Theorem~{\rm\ref{thm-main0}}, the
second fundamental form $\Pi$ at a point $p\in M$ is uniquely
determined, up to a unitary transformation of the target space, by
the tangential pseudoconformal curvature tensor $S$ of $M$ at $p$}
(and, hence, does not depend on $f$).

We also show (see Corollary~\ref{thm-secform}(i)) that {\em if
$d<n$ and $M$ is the sphere, then the second fundamental form of
$f$ vanishes identically}. Combining the latter result with
Theorem~\ref{thm-main3} above ($s=d$ in this case), we recover the
result of Faran {\it et al} mentioned above, that any CR-immersion
of the sphere in $\bC^{n+1}$ into the sphere in $\bC^{n+d+1}$,
with $d<n$, is equivalent (after composing to the left and right
with automorphisms of the spheres) to the linear embedding.

In \S\ref{sec-induced} we express the forms defining the
pseudoconformal connection of $\hat M$ pulled back to $M$ in terms
of the corresponding forms for $M$ and the second fundamental form
$\Pi$ of $f$ under the additional assumption that $\Pi$ is
nondegenerate. If it is not, higher order covariant derivatives of
$\Pi$ are needed; this is dealt with in \S \ref{sec-degsff}. We
begin the latter section by showing how covariant derivatives of
$\Pi$ can be used to determine the spaces $E_k$ defined in
\eqref{eq-ksdeg1} (see Proposition~\ref{pro-basic}). We then prove
(Theorem~\ref{thm-conn2}) that, if $f$ is $k_0$-nondegenerate,
then the pseudoconformal connection of $\hat M$ pulled back to $M$
is uniquely determined by the covariant derivatives of $\Pi$ up to
order $k_0-1$. The second main result in \S \ref{sec-degsff}
(Theorem~\ref{thm-secform2}) then states that equality of the
latter derivatives for two immersions always holds, after possibly
a unitary change of adapted coframes, provided that the
codimension $d<n/2$. An important technical point here is to
obtain a commuting relation between the covariant derivatives
(Lemma~\ref{lem-covdercom}). Finally, in \S\ref{sec-adapt} we
recall from \cite{CM} how adapted $Q$-frames on a sphere and their
Maurer-Cartan forms are related to the pseudoconformal connection
forms (see \eqref{eq-piphi}). We then complete the proof of
Theorem~\ref{thm-main3} in the case $s=0$, i.e.\ when $f$ is
finitely nondegenerate at some point. In this case the results of
\S\ref{sec-degsff} yield that the covariant derivatives of $\Pi$
and hence the pulled back pseudoconformal connection does not
depend on $f$. From \eqref{eq-piphi} we further conclude that also
Maurer-Cartan forms for associated adapted $Q$-frames do not
depend on $f$ and the proof is completed as in \cite{We79} by
general ODE arguments. \S\ref{sec-deg} is mainly devoted to the
proof of Theorem~\ref{thm-main2} which is obtained as a
consequence of the fact that the span of certain vectors in a
suitable adapted $Q$-frame is independent of the reference point.
Theorem~\ref{thm-main2} is then used to complete the proof of
Theorem~\ref{thm-main3} by reducing the general case to the case
$s=0$ treated in \S\ref{sec-adapt}.

\section{Preliminaries}\Label{sec-prelim}

Let $M$ be a strictly pseudoconvex CR-manifold (which in this
paper will always be understood to be of hypersurface type) of
dimension $2n+1$. We shall write
$$T^cM:=TM\cap iTM\subset TM, \quad \mathcal V=T^{0,1}M:=\{X+iJX : X\in T^cM \} \subset \C TM:=TM\otimes\C$$
for its {\em maximal complex tangent bundle} and {\em CR-bundle} respectively
which are both complex rank $n$ bundles. Here $J\colon T^cM\to T^cM$ is the complex structure.
We also consider the cotangent bundles
\begin{equation}\Label{eq-bundles}
T^0M:=(\mathcal V\oplus\bar {\mathcal V})^\perp,\ T'M:=\mathcal
V^\perp. \end{equation} Thus, $T^0M$ and $T'M$ are rank one and
rank $n+1$ subbundles of $\bC T^* M$ respectively with $T^0\subset T'M$. The
bundle $T'M$ is called the {\em holomorphic} or $(1,0)$ {\em cotangent bundle} of $M$.
As usual, a section of $\Lambda^p(T'M)$ is called a {\em $(p,0)$-form} on $M$. A real
nonvanishing section $\theta$ of $T^0M$ is called a {\it contact form}.
A choice of a contact form defines uniquely a real vector field
$T$, the {\it characteristic} (or {\it Reeb}) vector field of
$\theta$ (cf. e.g.\ \cite{Ho98}), by
\begin{equation}\Label{eq-charvf}
T\lrcorner d\theta=0,\ \left<\theta,T\right>=1,\end{equation}
where $\lrcorner$ denotes contraction (or interior
multiplication).
Indeed, since $d\theta$ is a degenerate $2$-form on $TM$ but
nondegenerate on the real hyperplanes defined by $\theta=0$ in $TM$,
one can always find $T$ satisfying \eqref{eq-charvf} in the kernel of $d\theta$.

We will follow the notation of \cite{CM} and \cite{We78},
in particular, we use the summation convention
and small Greek indices will always run over the set $\{1,\ldots, n\}$.
A typical tensor will be written as $S_\a{}^{\b}{}_{\mu}{}_{\bar\nu}$,
where an index without (resp. with) conjugation indicates $\C$-linear
(resp. $\C$-antilinear) dependence in the corresponding argument.
Here such a tensor $S_\a{}^{\b}{}_{\mu}{}_{\bar\nu}$ can be considered as
an $\R$-multilinear complex-valued function on $\3V\times\3V^*\times\3V\times\3V$.
The tensors will not be necessarily symmetric in their indices
and hence the order of the indices will be important and will be explicitly
indicated.
Simultaneous conjugation of all indices corresponds to conjugation of the
tensor:
$S_{\bar\a}{}^{\bar\b}{}_{\bar\mu}{}_{\nu}=
\1{S_\a{}^{\b}{}_{\mu}{}_{\bar\nu}}$.
On the other hand, there is no a priori relation e.g.\ between
$S_\a{}^{\b}{}_{\mu}{}_{\bar\nu}$ and $S_\a{}^{\b}{}_{\mu}{}_{\nu}$.
The Levi form matrix $(g_{\a\bar\b})$ of $M$ (relative to a given contact form) and its inverse
$(g^{\a\bar\b})$ will be used to raise and lower indices (without changing their order):
$S_\a{}_{\bar\b}{}_{\mu}{}_{\bar\nu}=g_{\gamma\bar\beta} S_\a{}^{\gamma}{}_{\mu}{}_{\bar\nu}$.
More generally, the same notation will be used for indexed functions on $M$
that may not necessarily transform as tensors, e.g.\ for connection matrices etc.

If we choose a basis $L_\a$, $\a=1,\ldots, n$, of $(1,0)$ vector fields (i.e.\ sections of
$T^{1,0}M=\1{\3V}$), so that $(T,L_\a,L_{\bar \a})$ is a
frame for $\bC TM$, then the first equation in (\ref{eq-charvf}) is equivalent to
\begin{equation}\Label{eq-dtheta}
d\theta=ig_{\a\bar\b}\theta^\a\wedge\theta^{\bar\b},
\end{equation}
$(g_{\a\bar\b})$ is the (Hermitian) Levi form matrix as above
and $(\theta,\theta^\a,\theta^{\bar\a})$ is the coframe
(i.e.\ a collection of linearly independent $1$-forms spanning $\C T^*M$)
dual to $(T,L_\a,L_{\bar\a})$ (for brevity, we shall say
that $(\theta,\theta^\a)$ is the coframe dual to  $(T,L_\a)$).
Note that $\theta$ and $T$ are real whereas $\theta^\a$ and $L_\a$
always have nontrivial real and imaginary parts.

Following Webster~\cite{We78}, we call
a coframe $(\theta,\theta^\a)$ (and its dual frame $(T,L_\a)$), where
$\theta$ is a contact form, {\em admissible} if
\eqref{eq-dtheta} holds or, equivalently, if $T$
is characteristic for $\theta$ in the sense of \eqref{eq-charvf}.
Observe that (by the uniqueness of the Reeb vector field)
for a given contact form $\theta$ on $M$, the 
admissible coframes are determined up to transformations
\begin{equation}\Label{eq-adcof}
\tilde \theta^\a= {u_\b}^\a \theta^\b,
\quad ({u_\b}^\a)\in \GL(\bC^n).
\end{equation}
Every choice of a contact form $\theta$ on $M$ is called {\em pseudohermitian structure}
and defines a Hermitian metric on ${\mathcal V}$ (and on $\1{\3V}$) via the (positive-definite) Levi form.
For every such $\theta$, Tanaka~\cite{Ta75} and Webster~\cite{We78}
defined a {\em pseudohermitian connection} $\nabla$ on $\bar {\mathcal V}$ (and also on $\C TM$)
which is expressed relative to an admissible coframe $(\theta,\theta^\a)$ by
\begin{equation}\Label{eq-con}
\nabla L_\a:={\o_\a}^\b\otimes L_\b,
\end{equation}
where the $1$-forms ${\o_\a}^\b$ on $M$ are uniquely determined
by the conditions
\begin{equation}\Label{eq-consymmetry}
\begin{aligned}
d\theta^\b&=\theta^\a\wedge {\o_\a}^\b \mod \theta\wedge\theta^{\bar\a},\\
dg_{\a\bar\b}&=\o_{\a\bar\b}+ \o_{\bar\b\a}.
\end{aligned}
\end{equation}
The first condition in \eqref{eq-consymmetry} can be rewritten as
\begin{equation}\Label{eq-conmatrix}
d\theta^\b=\theta^\a\wedge {\o_\a}^\b + \theta\wedge\tau^\b,
\quad \tau^\b = A^\b{}_{\bar\nu}\theta^{\bar\nu},
\quad A^{\a\b} = A^{\b\a}
\end{equation}
for suitable uniquely determined torsion matrix $(A^\b{}_{\bar\a})$,
where the last symmetry relation holds automatically (see \cite{We78}).
(More precisely, the forms ${\o_\a}^\b$ and $\tau^\b$ are first defined in \cite{We78}
on the principal bundle $P$ of all admissible coframes $(\theta,\theta^\a)$
on $M$ with fixed $\theta$ and are then pulled back to $M$
via a section of $P$ corresponding to a choice of such a coframe.)

The curvature of the pseudohermitian connection is given, in view of \cite[(1.27), (1.41)]{We78}, by
\begin{equation}\Label{eq-curv}
d{\o_\a}^\b-{\o_\a}^\gamma\wedge{\omega_\gamma}^\b=
{{R_{\a}}^\b}_{\mu\bar\nu}
\theta^\mu\wedge\theta^{\bar\nu}
+W_\a{}^\b{}_\mu\theta^\mu\wedge\theta
- W^\b{}_{\a\bar\nu}\theta^{\bar \nu}\wedge \theta
+i \theta_\a\wedge\tau^{\b}
- i \tau_\a\wedge\theta^{\b},
\end{equation}
where the functions ${{R_{\a}}^\b}_{\nu\bar\mu}$ and $W_\a{}^\b{}_\nu$
represent the pseudohermitian curvature of $(M,\theta)$.
It has been noticed by Lee \cite{Lee1} that the components $W_\a{}^\b{}_\mu$
can be in fact obtained as covariant derivatives of the torsion matrix $A^\b{}_{\bar\a}$ in
\eqref{eq-conmatrix}. Here we denote the covariant differentiation operator
with respect to the pseudohermitian connection $\nabla$ also by $\nabla$
and its components by indices preceded by a semicolumn, where the index $0$ is used to denote the covariant
derivative with respect to $T$; thus, e.g.\
\begin{equation}\Label{eq-cov}
\nabla A^{\b}{}_{\bar\a}
= dA^{\b}{}_{\bar\a} +  A^\mu{}_{\bar\a} {\o_\mu}^\b  -  A^\b{}_{\bar\nu}   {\o_{\bar\a}}^{\bar\nu}
=  A^\b{}_{\bar\a;0}\theta + A^\b{}_{\bar\a;\nu}\theta^\nu + A^\b{}_{\bar\a;\bar\nu}\theta^{\bar\nu}.
\end{equation}
In this notation the above mentioned relation reads \cite[(2.4)]{Lee1}:
\begin{equation}\Label{eq-goodframe}
W_\a{}^\b{}_\mu=A_{\a\mu;}{}^\b,
\quad W^\b{}_{\a\bar\nu}=A^\b{}_{\bar\nu}{}_{;\a}.
\end{equation}

We shall also need the Chern-Moser coframe bundle $Y$ over $M$.
Recall \cite[\S4]{CM} that $Y$ is the bundle of the coframes $(\o,\o^\a,\o^{\bar\a},\phi)$
on the real line bundle $\pi_E\colon E\to M$ (of all contact forms) satisfying
$d\o=ig_{\a\bar\b}\o^\a\wedge\o^{\bar\b}+\o\wedge\phi$,
where $\o^\a$ is in $\pi_E^* (T'M)$
and $\o$ is the canonical form on $E$ given by
$\o(\theta)(X):=\theta((\pi_E)_*X)$ for $\theta\in E$, $X\in T_\theta E$.
Similarly, canonical forms $\o,\o^\a,\o^{\bar\a},\phi$ are defined on $Y$
(here the same letters are used as for the coframe by a slight abuse of notation).
Chern and Moser \cite{CM} showed that
these forms can be completed to a natural parallelism on $Y$
given by the coframe of $1$-forms
\begin{equation}\Label{eq-CRpar}
(\omega,\omega^\a,\omega^{\bar\a},\phi,\phi_{\b}{}^{\a},\phi^{\a},
\phi^{\bar\a},\psi)
\end{equation}
defining the {\em pseudoconformal connection} on $Y$ and
satisfying the structure equations
(see \cite{CM} and its appendix)
\begin{equation}\Label{eq-cmw}
\begin{aligned}
& \phi_{\a\bar\b}+\phi_{\bar\b\a}=g_{\a\bar\b}\phi,\\
& d\omega=i\omega^\mu\wedge\omega_{\mu}+\omega\wedge\phi,\\
& d\omega^\a=\omega^{\mu}\wedge\phi_{\mu}{}^{\a}+\omega\wedge\phi^\a,\\
& d\phi=i\omega_{\bar\nu}\wedge\phi^{\bar\nu}
 +i\phi_{\bar\nu}\wedge\omega^{\bar\nu}+\omega\wedge\psi,\\
& d\phi_{\b}{}^\a=\phi_{\b}{}^\mu\wedge\phi_{\mu}{}^{\a}
 +i\omega_\b \wedge\phi^\a-i\phi_\b\wedge\omega^\a
 -i\delta_\b{}^\a\phi_\mu\wedge\omega^\mu
 -\frac{\delta_\b{}^\a}{2}\psi\wedge\omega+\Phi_\b{}^\a,\\
& d\phi^\a=\phi\wedge\phi^\a+\phi^\mu\wedge\phi_\mu{}^\a
 -\frac{1}{2}\psi\wedge\omega^\a+\Phi^\a,\\
& d\psi=\phi\wedge\psi+2i\phi^\mu\wedge\phi_\mu+\Psi,
\end{aligned}
\end{equation}
where the curvature $2$-forms $\Phi_\b{}^\a, \Phi^\a, \Psi$ can be decomposed as
\begin{equation}\Label{eq-CRtor}
\begin{aligned}
&\Phi_\b{}^\a=S_\b{}^\a{}_{\mu\bar\nu}\omega^\mu\wedge\omega^{\bar\nu}+
V_\b{}^{\a}{}_{\mu}\omega^\mu\wedge\omega+V^\a{}_{\b\bar\nu}\omega\wedge
\omega^{\bar\nu},\\
&\Phi^\a=V^\a{}_{\mu\bar\nu}\omega^\mu\wedge\omega^{\bar\nu}
+P_\mu{}^\a\omega^\mu\wedge\omega
+Q_{\bar\nu}{}^\a\omega^{\bar\nu}\wedge\omega,\\
&\Psi=-2iP_{\mu\bar\nu}\omega^\mu\wedge\omega^{\bar\nu}+R_\mu\omega^\mu\wedge\omega
+R_{\bar\nu}\omega^{\bar\nu}\wedge\omega,
\end{aligned}
\end{equation}
where the functions $S_\b{}^\a{}_{\mu\bar\nu}$, $V_\b{}^{\a}{}_{\mu}$,
$P_\mu{}^\a$, $Q_{\bar\nu}{}^\a$, $R_\mu$ together
represent the {\em pseudoconformal curvature} of $M$
(the indices of $S_\b{}^\a{}_{\mu\bar\nu}$ here are interchanged comparing to \cite{CM}
to make them consistent with indices of $R_\b{}^\a{}_{\mu\bar\nu}$ in \eqref{eq-curv}).
As in \cite{CM} we restrict our attention here to coframes
$(\theta,\theta^\a)$ for which the Levi form
$(g_{\a\bar\b})$ is constant. The 1-forms $\phi^{\a},\phi^{\bar\a},\phi_{\b}{}^{\a},\psi$
are uniquely determined by requiring the coefficients in
(\ref{eq-CRtor}) to satisfy certain symmetry and trace conditions
(see \cite{CM} and the appendix), e.g.
\begin{equation}\Label{eq-Strace}
S_{\a\bar\b\mu\bar\nu} = S_{\mu\bar\b\a\bar\nu}=S_{\mu\bar\nu\a\bar\b}=S_{\bar\nu\mu\bar\b\a},
\quad S_{\mu}{}^\mu{}_{\a\bar\b}=V_{\a}{}^\mu{}_\mu=P_\mu{}^\mu=0.
\end{equation}

Let us fix a contact form $\theta$ that defines a section $M\to E$.
Then any admissible coframe $(\theta,\theta^\a)$ for $T'M$
defines a unique section $M\to Y$ for which the pullbacks
of $(\o,\o^\a)$ coincide with $(\theta,\theta^\a)$
and the pullback of $\phi$ vanishes.
As in \cite{We78}, we use this section to pull back the forms (\ref{eq-CRpar}) to $M$.
We shall use the same notation for the pulled back forms on $M$
(that now depend on the choice of the admissible coframe).
With this convention, we have
\begin{equation}\Label{pbvanish}
\theta=\omega, \quad \theta^\a=\omega^\a, \quad \phi=0
\end{equation}
on $M$.
Now, in view of \cite[(3.8)]{We78}, the pulled back tangential pseudoconformal curvature tensor
$S_{\a}{}^{\b}{}_{\mu\bar\nu}$ can be obtained
from the tangential pseudohermitian curvature tensor ${{R_\a}^\b}_{\mu\bar\nu}$ in \eqref{eq-curv} by
\begin{equation}\Label{eq-curvid}
S_{\a\bar\b\mu\bar\nu}=R_{\a\bar\b\mu\bar\nu}
-\frac{R_{\a\bar\b}g_{\mu\bar\nu}
+R_{\mu\bar\b} g_{\a\bar\nu}
+R_{\a\bar\nu}g_{\mu\bar\b}
+R_{\mu\bar\nu}g_{\a\bar\b}}{n+2}
+\frac{R
(g_{\a\bar\b} g_{\mu\bar\nu}+
g_{\a\bar\nu} g_{\mu\bar\b})}
{(n+1)(n+2)},
\end{equation}
where
\begin{equation}\Label{Ricci}
R_{\a\bar\b}:=R_\mu{}^\mu{}_{\a\bar\b} \text{ and } R:=R_{\mu}{}^\mu
\end{equation}
are respectively the {\em pseudohermitian Ricci} and {\em scalar
curvature} of $(M,\theta)$. Formula \eqref{eq-curvid} expresses
the fact that $S_{\a\bar\b\mu\bar\nu}$ is the ``traceless
component'' of $R_{\a\bar\b\mu\bar\nu}$ with respect to the
natural decomposition of the space of all tensors
$T_{\a\bar\b\mu\bar\nu}$ with the symmetry condition as for
$S_{\a\bar\b\mu\bar\nu}$ in \eqref{eq-Strace} into the direct sum
of the subspace of such tensors of trace zero (i.e.\
$T_{\mu}{}^\mu{}_{\a\bar\b}=0$) and the subspace of ``multiples of
the Levi form'', i.e.\ tensors of the form
\begin{equation}\Label{mult}
T_{\a\bar\b\mu\bar\nu}
= H_{\a\bar\b}g_{\mu\bar\nu}+H_{\mu\bar\b} g_{\a\bar\nu}
+H_{\a\bar\nu}g_{\mu\bar\b}+H_{\mu\bar\nu}g_{\a\bar\b},
\end{equation}
where $(H_{\a\bar\b})$ is any Hermitian matrix.
We shall call two tensors as above {\em conformally equivalent}
if their difference is of the form \eqref{mult}.
In this terminology, the right hand side of \eqref{eq-curvid}
(together with \eqref{Ricci})
gives for any tensor $R_{\a\bar\b\mu\bar\nu}$
(with the above symmetry relations) its {\em traceless component}
which is the unique tensor of trace zero
that is conformally equivalent to $R_{\a\bar\b\mu\bar\nu}$.

The following result establishes relations between
pseudoconformal and pseudohermitian connection forms
and is alluded to in \cite{We78}.

\begin{Pro}\Label{prop1}
Let $M$ be a strictly pseudoconvex
CR-manifold of hypersurface type of CR-dimension $n$,
let $\omega_\b{}^\a$, $\tau^\a$ be defined by
{\rm(\ref{eq-consymmetry}--\ref{eq-conmatrix})}
with respect to an admissible coframe $(\theta,\theta^\a)$
and let $\phi_\b{}^\a$, $\phi^\a$, $\psi$ be
the forms in {\rm\eqref{eq-CRpar}} pulled back
to $M$ using $(\theta,\theta^\a)$ as above.
Then we have the following relations:
\begin{equation}\Label{eq-CRherm1}
\phi_\b{}^\a=\omega_\b{}^\a+D_\b{}^\a\theta, \quad
\phi^\a=\tau^\a + D_\mu{}^\a\theta^\mu  + E^\a\theta, \quad
\psi=iE_\mu\theta^\mu-iE_{\bar\nu}\theta^{\bar\nu}+B\theta
\end{equation}
where
\begin{equation}\Label{eq-CRherm2}
\begin{aligned}
&D_{\a\bar\b}:=\frac{iR_{\a\bar\b}}{n+2}-\frac{iRg_{\a\bar\b}}{2(n+1)(n+2)},\\
&E^\a:=\frac{2i}{2n+1}\big(A^{\a\mu}{}_{;\mu}-D^{\bar\nu\a}{}_{;\bar\nu}\big),\\
&B:=\frac{1}{n}\big( E^\mu{}_{;\mu} + E^{\bar\nu}{}_{;\bar\nu} -2A^{\b\mu}A_{\b\mu}
+2D^{\bar\nu\a} D_{\bar\nu\a} \big).
\end{aligned}
\end{equation}
\end{Pro}

\begin{proof}
The formulas for $\phi_\b{}^\a$ and $D_{\a\bar\b}$ were proved
in \cite{We78}. The
formula for $\phi^\a$ follows from the third equation
in (\ref{eq-cmw}) and (\ref{eq-conmatrix}). Indeed, these two
equations yield
\begin{equation*}
\theta^\a\wedge {\o_\a}^\b +\theta\wedge\tau^\b=
\theta^\a\wedge\phi_\a{}^\b+\theta\wedge\phi^\b.
\end{equation*}
Substituting the formula for $\phi_\b{}^\a$ in \eqref{eq-CRherm1}, we obtain
\begin{equation}
\theta\wedge\tau^{\b} = D_\a{}^\b\theta^\a\wedge\theta + \theta\wedge\phi^\b,
\end{equation}
which implies the formula for $\phi^\a$ in \eqref{eq-CRherm1} with
some $E^\a$.
Similarly, the formula for $\psi$ in \eqref{eq-CRherm1} with some $B$ follows
from equating the coefficients of $\theta$ in the pulled back fourth equation of (\ref{eq-cmw})
and using \eqref{pbvanish} (whence $d\phi=0$ on $M$).

To obtain the formula for $E^\a$ in \eqref{eq-CRherm2}, we substitute the formulas \eqref{eq-CRherm1}
for $\phi_\b{}^\a,\phi^\a,\psi$ in the pulled back
sixth equation of (\ref{eq-cmw}) and use \eqref{eq-dtheta}, \eqref{eq-conmatrix},
the covariant derivative \eqref{eq-cov} (and the analogue for $D_\b{}^\a$)
and the formula for $\Phi^\a$ in \eqref{eq-CRtor}:
\begin{equation}\Label{eq-DA}
\nabla A^{\a}{}_{\bar\nu}\wedge
\theta^{\bar\nu}+ \nabla D_\b{}^\a\wedge\theta^\b+ig_{\mu\bar\nu}E^\a\theta^\mu
\wedge\theta^{\bar\nu}=-\frac{1}{2}\psi\wedge\theta^\a+V^\a{}_{\mu\bar\nu}
\theta^\mu \wedge\theta^{\bar\nu} \mod\theta.
\end{equation}
By identifying the coefficient in front of
$\theta^\mu\wedge\theta^{\bar\nu}$ in \eqref{eq-DA} and using the formula for $\psi$
in \eqref{eq-CRherm1}, we obtain
\begin{equation*}
A^{\a}{}_{\bar\nu;\mu}-D_\mu{}^\a{}_{;\bar\nu}+ig_{\mu\bar\nu}E^\a=
-\frac{1}{2}iE_{\bar\nu}\delta_{\mu}{}^\a+V^\a{}_{\mu\bar\nu}.
\end{equation*}
The formula for $E^\a$ in \eqref{eq-CRherm2} is now obtained by
summing over $\mu$ and $\bar\nu$ and using the trace condition
$V^\a{}_\mu{}^\mu=0$.
Similarly, the formula for $B$ follows by
substituting the formula for $\psi$ in pulled back last equation of (\ref{eq-cmw}) (mod $\theta$)
and using the trace condition $P_\nu{}^\nu=0$.
\end{proof}

\section{Submanifolds of CR-manifolds; the second fundamental form}
\Label{sec-subs}

Let $M$ be a strictly pseudoconvex CR-manifold  of dimension $2n+1$ as before
and $f\colon M\to \hat M$ be an CR-immersion of $M$
into another  strictly pseudoconvex CR-manifold $\hat M$ of dimension $2{\hat n}+1$,
with rank ${\hat n}$ CR-bundle $\hat {\mathcal V}$.
Our arguments in the sequel will be of local nature
and hence we shall assume $f$ to be an embedding.
We shall use a $\hat{}$ to denote objects associated to $\hat M$.
Capital Latin indices $A, B,$ etc,
will run over the set $\{1,2,\ldots,{\hat n}\}$ whereas Greek
indices $\a,\b$, etc, will run over $\{1,2,\ldots, n\}$ as above.
Moreover, we shall let small Latin indices $a,b,$ etc, run over
the complementary set $\{n+1,n+2,\ldots, {\hat n}\}$.
Since $\hat M$ is strictly pseudoconvex and $f$ an embedding, it
is well known that for any contact form $\hat\theta$ on $\hat M$
the pullback $f^*(\hat\theta)$ (which for a CR-mapping $f$ is
always a section of $T^0M$) is nonvanishing and, hence, a contact
form on $M$ (In general, $f^*(\hat\theta)$ may vanish,
e.g. if $f(M)$ is contained in a complex-analytic subvariety of $\hat M$).
We shall always choose the coframe $(\hat\theta,\hat\theta^A)$ on $\hat M$
such that the pullback of $(\hat\theta,\hat\theta^\a)$ is a coframe for $M$
and hence drop the $\hat{}$ over the frames and coframes if there is no ambiguity.
It will be clear from the context if a form is pulled back to $M$ or not.

We shall identify $M$ with the submanifold $f(M)$ of $\hat M$
and write $M\subset\hat M$.
Then $\mathcal V$ becomes a rank $n$ subbundle of
$\hat {\mathcal V}$ along $M$. It follows that the (real) codimension of
$M$ in $\hat M$ is $2({\hat n}-n)$ and that there is a rank $({\hat n}-n)$
subbundle $N'M$ of $T'\hat M$ along $M$
consisting of $1$-forms on
$\hat M$ whose pullbacks to $M$ (under $f$) vanish.
We shall call $N'M$ the
{\em holomorphic conormal bundle of $M$} in $\hat M$.
We shall say that the pseudohermitian structure $(\hat M,\hat\theta)$
(or simply $\hat\theta$) is
{\em admissible for the pair $(M,\hat M)$} if the characteristic vector
field $\hat T$ of $\hat \theta$ is tangent to $M$ (and hence coincides
with the characteristic vector field of the pullback of $\hat\theta$).
This is equivalent, as the reader can easily verify, to requiring that for
any admissible coframe $(\hat \theta,\hat \theta^A)$ on $\hat M$,
where $A=1,\ldots, {\hat n}$, the holomorphic conormal bundle $N'M$
is spanned by suitable linear combinations of the $\hat\theta^A$.
It is easily seen that not all
contact forms $\hat\theta$ are admissible for $(M,\hat M)$.
However, we have the following statement:

\begin{Lem}
Let $M\subset \hat M$ be as above. Then any contact form on $M$
can be extended to a contact form $\theta$ in a neighborhood of $M$ in $\hat M$
such that $\theta$ is admissible for $(M,\hat M)$.
Moreover, the $1$-jet of $\theta$ is uniquely determined on $M$.
\end{Lem}

\begin{proof}
Let $\theta$ be any fixed extension of the given contact form on $M$ to a neighborhood of $M$ in $\hat M$.
Any other extension is clearly of the form $\tilde{\theta}=u\theta$,
where $u$ is a smooth function on $\hat M$ near $M$ with $u|_M \equiv 1$.
Let $T$ be the characteristic vector field of the restriction of $\theta$ to $M$.
Then $\tilde\theta$ is admissible for $(M,\hat M)$ if and only if $T\lrcorner d\tilde\theta=0$,
i.e.\ if $T\lrcorner d\theta - du = 0$ along $M$.
By the assumptions, the latter identity holds when pulled back to $M$.
Now it is clear that there exists unique choice of $du$ along $M$
for which it holds also in the normal direction.
The required function $u$ can be now constructed
in local coordinate charts and glued together via partition of unity.
The proof is complete.
\end{proof}

By taking admissible coframes as in \S\ref{sec-prelim} and using the Gram-Schmidt algorithm,
we obtain the following corollary,
where we take a little more care to distinguish between $M$ and its image $f(M)$ in $\hat M$.

\begin{Cor}\Label{thm-admada}
Let $M$ and $\hat M$ be strictly pseudoconvex CR-manifolds of dimensions $2n+1$ and
$2{\hat n}+1$ respectively and $f\colon M\to \hat M$ be a CR-embedding.
If $(\theta,\theta^\a)$ is any admissible coframe on $M$,
then in a neighborhood of any point $\hat p\in f(M)$ in $\hat M$
there exists an admissible coframe $(\hat \theta,\hat
\theta^A)$ on $\hat M$ with
$f^*(\hat\theta,\hat \theta^\a,\hat\theta^a)=(\theta,\theta^\a,0)$.
In particular, $\hat\theta$ is admissible for the pair $(f(M),\hat M)$,
i.e.\ the characteristic vector field $\hat T$ is tangent to $f(M)$.
If the Levi form of $M$ with respect to $(\theta,\theta^\a)$ is $\delta_{\a\bar\b}$,
then $(\hat \theta,\hat \theta^A)$ can be chosen such that the Levi form of $\hat M$
relative to it is also $\delta_{A\bar B}$.
With this additional property,
the coframe $(\hat \theta,\hat \theta^A)$
is uniquely determined along $M$ up to unitary transformations in $\U(n)\times\U(\hat n-n)$.
\end{Cor}

Let us fix an admissible coframe $(\theta,\theta^\a)$ on $M$
and let $(\hat \theta,\hat \theta^A)$ be an admissible coframe on
$\hat M$ near $f(M)$.
We shall say that
$(\hat \theta,\hat \theta^A)$ is {\em adapted} to
$(\theta,\theta^\a)$ on $M$ (or simply to
$M$ if the coframe on $M$ is understood) if it satisfies the
conclusion of Corollary~{\rm\ref{thm-admada}}
with the requirement there for the Levi form.

The fact that $(\theta,\theta^A)$ (where we omit a $\hat{}$ )
is adapted to $M$ implies, in view of the construction
(\ref{eq-consymmetry}), that if the
pseudohermitian connection matrix of $(\hat M,\theta)$ is $\hat\omega_B{}^A$,
then that of $(M,\theta)$ is (the pullback of) $\hat\omega_\b{}^\a$.
Similarly, the pulled back torsion matrix $\hat\tau^\a$ is $\tau^\a$.
Hence omitting a $\hat{}$ over these pullbacks will not cause any ambiguity
and we shall do it in the sequel.
By our normalization of the Levi form, the second equation in (\ref{eq-consymmetry}) reduces to
\begin{equation}\Label{eq-consym2}
\omega_{B\bar A}+\omega_{\bar AB}=0,
\end{equation}
where as before $\o_{\bar AB}=\1{\o_{A\bar B}}$.

The matrix of $1$-forms $({\omega_\a}^b)$ pulled back to $M$ defines the
{\it second fundamental form} of $M$ (or more precisely
of the embedding $f$). It was used e.g.\ in \cite{We79,Faran90}
in the
study of mapping problems. Since $\theta^b$ is 0 on $M$, we deduce
by using equation (\ref{eq-conmatrix}) that, on $M$,
\begin{equation}
{{\omega_\a}^b}\wedge \theta^\a
+ \tau^b\wedge\theta=0,
\end{equation} which
implies that
\begin{equation}\Label{eq-secform}
{\omega_\a}^b={{\omega_\a}^b}_\b\,\theta^\b,
\quad \omega_\a{}^b{}_\b=\omega_\b{}^b{}_\a,
\quad \tau^b=0.
\end{equation}

We now relate so defined intrinsic second fundamental form
$(\o_\a{}^b{}_\b)$
with the extrinsic one defined in \S\ref{sec-main}
in case $\hat M$ is embedded as a real hypersurface in $\C^{\hat n+1}$.
Given any admissible contact form $\theta$ for $(M,\hat M)$,
we can choose a defining function of $\hat M$ near a point $p=\hat p\in M$
such that $\theta=i\bar\partial\hat\rho$ on $\hat M$,
i.e. in local coordinates $Z'$ in $\C^{\hat n+1}$ vanishing at $p$ we have
\begin{equation}
\theta=i\sum_{k=1}^{{\hat n}+1}\frac{\partial\hat\rho}{\partial \bar Z'_k}\, d\bar Z'_k,
\end{equation}
where we pull back the forms $d \bar Z'_1,\ldots, d \bar Z'_{{\hat n}+1}$ to $\hat M$.
Given further a coframe $(\theta,\theta^A)$ on $\hat M$ near $p$ adapted to $M$
and its dual frame $(T,L_A)$, we have
\begin{equation}\Label{Lb}
L_\b (\hat\rho_{\bar Z'}\circ f) = -i L_\b \lrcorner d\theta =
g_{\b\bar C}\theta^{\bar C} = g_{\b\bar \gamma}\theta^{\bar \gamma}.
\end{equation}
Recall that we are in fact assuming that the Levi form has been normalized, i.e.\
$g_{A\bar B }=\delta_{A\bar B}$, even though we retain the
notation $g_{A\bar B }$.
Conjugating \eqref{Lb} we see that the subspace
$E_1(p)\subset T'_p\C^{\hat n+1}$ in \eqref{eq-ksdeg1} is spanned by
$(\theta,\theta^\a)$, where we use the standard identification
$T'_{p} \hat M\cong T'_{p}\bC^{{\hat n}+1}$.
Applying $L_\a$ to both sides of \eqref{Lb}
and using the analogue of (\ref{eq-consymmetry}) for $\hat M$ and
(\ref{eq-secform}), we conclude that
\begin{equation}\Label{eq-omega}
L_\a L_\b (\hat\rho_{\bar Z'}\circ f) =
g_{\b\bar \gamma} L_{\a}\lrcorner d\theta^{\bar \gamma}
= -\omega_{\bar a\b\a} \theta^{\bar a} = \omega_{\a \bar a\b} \theta^{\bar a}
\quad\mod \theta, \theta^{\bar\a},
\end{equation}
where we have used \eqref{eq-consym2} for the last identity.
Comparing with the extrinsic definition of the second fundamental form
(\ref{eq-defsec}) and identifying the spaces in \eqref{eq-defsec} and \eqref{SFF}
via the Levi form of $\hat M$ as explained in \S\ref{sec-main},
we conclude that
\begin{equation}\Label{eq-secid}
\Pi(L_\a,L_\b)=\omega_\a{}^a{}_\b \, L_a,
\end{equation}
where with have identified $L_a$ with its
equivalence classe in $T^{1,0}_{\hat p}\hat M/T^{1,0}_pM$.
Conjugating \eqref{eq-omega} and comparing with \eqref{eq-ksdeg1}
we see that the space $E_2=E_2(p)$ is spanned (via the identification above)
by the forms
\begin{equation}\Label{e2-forms}
\theta,\theta^\a,\omega_{\bar\a a \bar\b}\theta^a.
\end{equation}
To interpret the higher order spaces $E_k(p)$ in terms of the
second fundamental form, we shall need the covariant differentials
of the second fundamental form with respect to the induced
pseudoconformal connection, which will be introduced in section
\ref{sec-induced}. The discussion will therefore be postponed
until section \ref{sec-deg}.

As a byproduct of the relation \eqref{eq-secid} we see that the
bilinear map $T^{1,0}_pM\times T^{1,0}_pM\to T^{1,0}_p\hat M/T^{1,0}_pM$
defined by $(\o_\a{}^a{}_\b)$ is independent of the choice of the adapted coframe $(\theta,\theta^A)$
in case $\hat M$ is locally CR-embeddable (in $\C^{\hat n+1}$). More generally, if $\hat M$ is not CR-embeddable,
one can still use approximate embeddings constructed in \cite{KZ}
(or direct computation of the change of $\o_\a{}^a{}_\b$)
to obtain the same conclusion.

\section{CR analogue of the Gauss equation and applications}\Label{sec-gauss}
The classical Gauss equation in Riemannian geometry relates
the Riemannian curvature tensors of a manifold and its submanifold with the second fundamental form
composed with the Riemannian metric (see e.g.\ \cite{KN2}).
Our goal in this section will be to establish
the pseudohermitian and pseudoconformal analogues of the Gauss equation
and to apply them to obtain rigidity properties of the second fundamental form.

As before we fix a coframe $(\theta,\theta^A)$ adapted to $M$.
We first compare pseudohermitian curvature tensors $R_{\a}{}^\b{}_{\mu\bar\nu}$ and
$\hat R_{A}{}^B{}_{C\bar D}$ of $(\hat M, \theta)$ and $(M,\theta)$ respectively.
By comparing (\ref{eq-curv}) and the corresponding equation for
$\hat R_{A}{}^B{}_{C\bar D}$ pulled back to $M$
and using $\hat\omega_\a{}^\b=\omega_\a{}^\b$,
$\hat\tau^\a=\tau^\a$ and $\hat W_\a{}^\b{}_\mu=W_\a{}^\b{}_\mu$
as a consequence of \eqref{eq-goodframe},
we conclude that on $M$,
\begin{equation}\Label{eq-curvMhatM}
\hat R_{\a}{}^\b{}_{\mu\bar \nu }\,\theta^\mu\wedge\theta^{\bar
\nu}+ \omega_\a{}^a\wedge \omega_a{}^\b=R_{\a}{}^\b{}_{\mu\bar\nu}\,
\theta^\mu\wedge\theta^{\bar\nu}.
\end{equation}
By using the symmetry
(\ref{eq-consym2}), we conclude that, on $M$, we have
\begin{equation}\Label{eq-curvMhatM2}
\hat R_{\a}{}_{\bar\b}{}_{\mu\bar\nu}\,\theta^\mu\wedge
\theta^{\bar\nu}- \omega_\a{}^a\wedge \omega_{\bar\b a} =
R_{\a}{}_{\bar \b}{}_{\mu\bar\nu}\,
\theta^\mu\wedge\theta^{\bar\nu}.
\end{equation}
This can also be written, in view of (\ref{eq-secform}),
after equating the coefficients of $\theta^\mu\wedge\theta^{\bar\nu}$ as
\begin{equation}\Label{eq-relcurv}
\hat R_{\a\bar\b\mu\bar\nu}
= R_{\a\bar\b\mu\bar\nu}
+g_{a\bar b}\,\omega_\a{}^a{}_\mu\, \omega_{\bar\b}{}^{\bar b}{}_{\bar\nu} .
\end{equation}
The identity \eqref{eq-relcurv}
relates the tangential pseudohermitian curvature tensors of $M$ and $\hat M$
with the second fundamental form of the embedding $M$ into $\hat M$
and hence can be viewed as pseudohermitian analogue of the Gauss equation.
We state it in an invariant form using the previously established relation
\eqref{eq-secid} between the extrinsic and intrinsic second fundamental forms
$\Pi$ and $(\o_\a{}^a{}_\b)$ given respectively by (\ref{eq-defsec}--\ref{SFF})
and \eqref{eq-secform}. For this, we view both pseudohermitian curvature tensors
as $\R$-multilinear functions
$$R,\hat R\colon T^{1,0}M\times T^{1,0}M\times T^{1,0}M\times T^{1,0}M\to\C$$
that are $\C$-linear in the even and $\C$-antilinear in the odd numbered arguments.
We further identify the quotient space $T^{1,0}_p\hat M/T^{1,0}_p M$ for $p\in M$
with the orthogonal complement of $T^{1,0}_p M$ in $T^{1,0}_p\hat M$
with respect to the Levi form of $\hat M$ relative to $\theta$
and then use this Levi form to define a canonical Hermitian scalar product $\langle,\rangle$ on
$T^{1,0}_p\hat M/T^{1,0}_p M$. The identity \eqref{eq-relcurv} yields now the following statement.

\begin{Pro}[Pseudohermitian Gauss equation]\Label{psherm-gauss}
Let $M\subset \hat M$ be as above and $\theta$ be a contact form on $\hat M$
that is admissible for the pair $(M,\hat M)$. Then, for every $p\in M$, the following holds:
\begin{equation}\Label{psh-gauss}
\hat R (X,Y,Z,V) = R (X,Y,Z,V) + \big\langle \Pi(X,Z), \Pi(Y,V) \big\rangle,
\quad X,Y,Z,V \in T^{1,0}_pM.
\end{equation}
\end{Pro}

We next turn to a pseudoconformal analogue of the Gauss equation.
It follows immediately from \eqref{eq-relcurv} and \eqref{psh-gauss}
by taking traceless components of both sides
as discussed in \S\ref{sec-prelim}.
We write $[T_{\a\bar\b\mu\bar\nu}]$ for the traceless component
of a tensor $T_{\a\bar\b\mu\bar\nu}$
that can be computed by the formulas (\ref{eq-curvid}--\ref{Ricci})
(with $R_{\a\bar\b\mu\bar\nu}$ replaced by $T_{\a\bar\b\mu\bar\nu}$).
Then \eqref{eq-curvid} can be rewritten as
\begin{equation}\Label{SR}
S_{\a\bar\b\mu\bar\nu} = [R_{\a\bar\b\mu\bar\nu}].
\end{equation}
By taking traceless components of both sides in \eqref{eq-relcurv}
and using \eqref{SR} we now obtain
\begin{equation}\Label{pre-gauss}
[\hat R_{\a\bar\b\mu\bar\nu}]
= S_{\a\bar \b\mu\bar\nu}
+[g_{a\bar b}\,\omega_\a{}^a{}_\mu\, \omega_{\bar\b}{}^{\bar b}{}_{\bar\nu}].
\end{equation}
Note that, in contrast to \eqref{SR}, the left hand side of \eqref{pre-gauss}
may not be, in general, equal to $\hat S_{\a\bar\b\mu\bar\nu}$
which is the (restriction of the) traceless component of $\hat R_{A\bar B C\bar D}$
with respect to the indices running from $1$ to $\hat n$.
However, we claim that $[\hat R_{\a\bar\b\mu\bar\nu}]=[\hat S_{\a\bar\b\mu\bar\nu}]$.
Indeed, according to the decomposition into zero trace component and a multiple of
the Levi form $(g_{A\bar B})$ on $\hat M$, the tensors
$\hat R_{A\bar B C\bar D}$ and $\hat S_{A\bar B C\bar D}$ are conformally equivalent
with respect to the Levi form $(g_{A\bar B})$, i.e.\ their difference is of the
form analogous to \eqref{mult} with Greek indices replaced by capital Latin ones.
Since the Levi form of $\hat M$ restricts to that of $M$,
the restrictions $\hat R_{\a\bar\b\mu\bar\nu}$ and $\hat S_{\a\bar\b\mu\bar\nu}$
are conformally equivalent with respect to $(g_{\a\bar\b})$ and the claim follows.
Hence \eqref{pre-gauss} immediately yields the desired relation
between the tangential pseudoconformal curvature tensors of $M$ and $\hat M$
and the second fundamental form:
\begin{equation}\Label{gauss}
[\hat S_{\a\bar\b\mu\bar\nu}]
= S_{\a\bar \b\mu\bar\nu}
+[g_{a\bar b}\,\omega_\a{}^a{}_\mu\, \omega_{\bar\b}{}^{\bar b}{}_{\bar\nu}],
\end{equation}
or,  using formulas (\ref{eq-curvid}--\ref{Ricci}) for the traceless part,
\begin{equation}\Label{gauss-explicit}
\begin{aligned}
S_{\alpha\bar\beta\mu\bar\nu}= &
\, \hat S_{\alpha\bar\beta\mu\bar\nu}
-\frac{
\hat S_\gamma{}^\gamma{}_{\alpha\bar\beta}\,g_{\mu\bar \nu}
+\hat S_\gamma{}^\gamma{}_{\mu\bar\beta}\,g_{\alpha\bar \nu}
+\hat S_\gamma{}^\gamma{}_{\alpha\bar\nu}\,g_{\mu\bar \beta}
+\hat S_\gamma{}^\gamma{}_{\mu\bar\nu}\,g_{\alpha\bar\beta}
}{n+2}\\
+&\frac{
\hat S_\gamma{}^\gamma{}_\delta{}^\delta \,
(g_{\alpha\bar\beta} g_{\mu\bar\nu}+ g_{\alpha\bar\nu} g_{\mu\bar \beta})
}{(n+1)(n+2)}
-g_{a\bar b}\omega_\alpha{}^a{}_\mu\,\omega_{\bar\b}{}^{\bar b}{}_{\bar \nu}\\
+&\frac{
\omega_\gamma{}^a{}_\alpha\,\omega^{\gamma}{}_{a \bar \beta}\, g_{\mu\bar \nu}
+\omega_\gamma{}^a{}_\mu\,\omega^{\gamma}{}_{a \bar \beta}\, g_{\alpha\bar \nu}
+\omega_\gamma{}^a{}_\alpha\,\omega^{\gamma}{}_{a \bar \nu}\, g_{\mu\bar \beta}
+\omega_\gamma{}^a{}_\mu\,\omega^{\gamma}{}_{a \bar \nu}\, g_{\alpha\bar\beta}
}{n+2}\\
-&\frac{
\omega_\gamma{}^a{}_\delta\,\omega^{\gamma}{}_a{}^{\delta}
(g_{\alpha\bar\beta} g_{\mu\bar\nu}+ g_{\alpha\bar\nu} g_{\mu\bar \beta})
}{(n+1)(n+2)}.
\end{aligned}
\end{equation}
When ${\hat n}=n+1$ and $\hat M$ is a sphere
(so that $\hat S_{\a\bar\b\mu\bar\nu} \equiv 0$),
the identity (\ref{gauss-explicit}) reduces to that
obtained by Webster \cite[(2.14)]{We79}.
As for the pseudohermitian curvature,
we also view the pseudoconformal curvature tensors
(as well as their zero trace components)
as $\R$-multilinear functions
$$S,\hat S\colon T^{1,0}M\times T^{1,0}M\times T^{1,0}M\times T^{1,0}M\to\C$$
but now they are independent of $\theta$.
Then, with the above notation, \eqref{gauss} yields
the following statement.

\begin{Pro}[Pseudoconformal Gauss equation]\Label{psconf-gauss}
For $M\subset \hat M$ as above and every $p\in M$, the following holds:
\begin{equation}\Label{psc-gauss}
\big[\hat S (X,Y,Z,V)\big]
= S (X,Y,Z,V) + \big[\big\langle \Pi(X,Z), \Pi(Y,V) \big\rangle\big],
\quad X,Y,Z,V \in T^{1,0}_pM.
\end{equation}
\end{Pro}

In the rest of this section we apply the Gauss equations (\ref{gauss}--\ref{psc-gauss})
to the case when $\hat M$ is the sphere and hence $\hat S_{A\bar B C\bar D} \equiv 0$.
Based on a lemma due to the second author \cite{Hu1}, we show
that the second fundamental form in this case
is uniquely determined up to unitary transformations
by the tangential pseudoconformal curvature $S_{\a\bar \b\mu\bar\nu}$
under suitable restrictions on the dimensions $n$ and $\hat n$.
We begin with an algebraic lemma.

\begin{Lem}\Label{alg}
Let $\Pi\colon \C^n\times\C^n \to \C^{\hat n-n}$ be a $\C$-bilinear map
and $\langle\cdot,\cdot \rangle$ an Hermitian scalar product on $\C^{\hat n-n}$.
Denote by $S:=\big[\big\langle \Pi(\cdot,\cdot), \Pi(\cdot,\cdot) \big\rangle\big]$
the traceless component with respect to a fixed Hermitian scalar product $(\cdot,\cdot)$ on $\C^n$ as above.
\begin{enumerate}
\item[(i)] If $\Pi\equiv 0$, then also $S\equiv 0$. The converse also holds under the assumption ${\hat n}-n<n$.
\item[(ii)] If ${\hat n}-n<n/2$, then $\Pi$ is uniquely determined by $S$ up to unitary transformations
of $\C^{\hat n-n}$, i.e.\ if $\tilde\Pi$ is another form, then
$S=\big[\big\langle \tilde\Pi(\cdot,\cdot), \tilde\Pi(\cdot,\cdot) \big\rangle\big]$
if and only if $\tilde\Pi=U\circ \Pi$ for some $U\in\GL(\hat n- n)$ preserving $\langle\cdot,\cdot \rangle$.
\end{enumerate}
\end{Lem}

\begin{proof}
The first statement of (i) is trivial.
By definition of the traceless component, $S\equiv 0$ implies
$\langle \Pi(X,X), \Pi(X,X) \big\rangle = (X,X)\, H(X,X)$ for
some Hermitian form $H(\cdot,\cdot)$ and all $X\in\C^n$.
Then the second statement of (i) follows directly from \cite[Lemma 3.2]{Hu1}.
To show (ii) observe that, for $\Pi$ and $\tilde\Pi$ as in the lemma, one similarly has
$$\langle \Pi(X,X), \Pi(X,X) \big\rangle - \big\langle \tilde\Pi(X,X), \tilde\Pi(X,X) \big\rangle = (X,X)\, H(X,X).$$
The restriction ${\hat n}-n<n/2$ permits now to conclude from \cite[Lemma 3.2]{Hu1} that $H(X,X)\equiv 0$.
By polarizing the left hand side and using the symmetries, we obtain
\begin{equation}\Label{eq-langle}
\langle \Pi(X,Y), \Pi(Z,V) \big\rangle = \big\langle \tilde\Pi(X,Y), \tilde\Pi(Z,V) \big\rangle,
\quad X,Y,Z,V \in \C^n,
\end{equation}
expressing the fact that the collections of vectors $\big(\Pi(X,Y)\big)_{X,Y}$
and $\big(\tilde \Pi(X,Y)\big)_{X,Y}$ have the same scalar products.
Hence they can be transformed into each other by a transformation $U\in \GL(\hat n- n)$
preserving $\langle\cdot,\cdot\rangle$.
Any such $U$ satisfies the required conclusion proving (ii).
\end{proof}

Example~\ref{ex-whitney} shows that (i) may not hold in case of equality $\hat n-n=n$.
Furthermore, the conclusion of (ii) may not hold under the weaker inequality $\hat n-n<n$
that was enough to obtain (i):

\begin{Ex}
For $n:=2$, $\hat n:=3$, consider the standard scalar products $(\cdot,\cdot)$
on $\C^2$ and $\langle\cdot,\cdot\rangle$ on $\C$
and define $\Pi$ and $\tilde\Pi$ by their associated quadratic forms
$Q:=4z_1^2 + z_1 z_2 + 4z_2^2$ and $\tilde Q:=4z_1^2 - z_1 z_2 + 4z_2^2$ respectively.
Then the identity
$$|4z_1^2 + z_1 z_2 + 4z_2^2|^2 - |4z_1^2 - z_1 z_2 + 4z_2^2|^2
 = 8(|z_1|^2 + |z_2|^2) (z_1 \bar z_2 + \bar z_1 z_2)$$
implies that
$\big[\big\langle \Pi(\cdot,\cdot), \Pi(\cdot,\cdot) \big\rangle\big]
=\big[\big\langle \tilde\Pi(\cdot,\cdot), \tilde\Pi(\cdot,\cdot) \big\rangle\big]$.
However, it is clear that the conclusion of (ii) does not hold for $\Pi$ and $\tilde\Pi$.
\end{Ex}

The transformation $U$ in Lemma~\ref{alg} is clearly uniquely determined on the image of $\Pi$.
We apply Lemma~\ref{alg} to compare second fundamental forms $\Pi=\Pi_p$ and $\tilde\Pi=\tilde\Pi_p$
at a point $p\in M$ for two given embedding $f,\tilde f\colon M\to \hat M$.
Then, if the image of $\Pi$ is of constant dimension,
in which case we say that $\Pi$ is {\em of constant rank},
the transformation $U$ can be chosen smoothly depending on $p$ (near any given point).
We then obtain, as a consequence of (\ref{psc-gauss}--\ref{gauss}) and Lemma~\ref{alg},
the following result.

\begin{Cor}\Label{thm-secform}
Let $f\colon M\to \hat M$ be a CR-embedding of
a strictly pseudoconvex CR-manifold of dimension $2n+1$
into the unit sphere $\hat M=\5S$ in $\C^{{\hat n}+1}$.
Denote by $(\omega_\a{}^a{}_\b)$
the second fundamental form matrix of $f$ relative
to an admissible coframe $(\theta,\theta^A)$ on $\hat M$ adapted to $f(M)$.
\begin{enumerate}
\item[(i)] If $\omega_\a{}^a{}_\b\equiv 0$, then  $M$ is locally CR-equivalent to the unit sphere in $\bC^{n+1}$.
The converse also holds under the assumption ${\hat n}-n<n$.
\item[(ii)] If ${\hat n}-n<n/2$, then, for any $p\in M$,
$(\omega_\a{}^a{}_\b)(p)$ is uniquely determined by
$(\theta,\theta^\a)$ up to unitary transformations of $\C^{\hat
n-n}$, i.e.\ for any other CR-embedding $\tilde f\colon M\to \hat
M$ and any coframe $(\tilde\theta,\tilde \theta^A)$ on $\hat M$
adapted to $\tilde f(M)$ with $\tilde
f^*(\tilde\theta,\tilde\theta^\a) = f^*(\theta,\theta^\a)$, one
has $\tilde\omega_\a{}^a{}_\b(p)=\omega_\a{}^a{}_\b(p)$ after
possibly a unitary change of $(\tilde\theta^a)$ near $p$.
Moreover, if $(\omega_\a{}^a{}_\b)$ is of constant rank near $p$,
there there is such a unitary change of $(\tilde\theta^a)$ near
$p$ for which $\tilde\omega_\a{}^a{}_\b=\omega_\a{}^a{}_\b$ near
$p$.
\end{enumerate}
\end{Cor}

\section{The induced pseudoconformal connection}\Label{sec-induced}
Let now $(\theta,\theta^A)$ be an adapted coframe for the pair
$(M,\hat M)$ as above and let
(\ref{eq-CRpar}) be the $1$-forms defining the parallelism on
the bundle $Y$ over $M$ pulled back to $M$ by the coframe
$(\theta,\theta^\a)$ as described in \S\ref{sec-prelim}.
We use $\hat{ }\,$ for the corresponding forms on $\hat M$
pulled back further to $M$,
where the indices run from $1$ to $\hat n$.
Recall that $(\o,\o^\a,\o^{\bar\a})=(\hat\o,\hat\o^\a,\hat\o^{\bar\a})=(\theta,\theta^\a,\theta^{\bar\a})$
and $\hat\o^a=0$ on $M$. In
view of Proposition \ref{prop1}, we do not expect
$(\phi_{\b}{}^{\a},\phi^{\a}, \psi)$ and
$(\hat \phi_{\b}{}^{\a},\hat \phi^{\a}, \hat \psi)$ to be
equal. However, since $\hat\omega_\b{}^\a=\omega_\b{}^\a$
and $\tau^\a=\hat\tau^\a$ (see \S\ref{sec-subs}),
Proposition \ref{prop1} implies
\begin{equation}\Label{eq-CRMhatM}
\hat \phi_\b{}^\a=\phi_\b{}^\a+C_\b{}^\a\theta,
\quad \hat \phi^\a=\phi^\a +C_\mu{}^\a\theta^\mu +F^\a\theta,
\quad \hat \psi=\psi+iF_\mu\theta^\mu-iF_{\bar\nu}\theta^{\bar\nu}+A\theta,
\end{equation}
where
\begin{equation*}
C_\b{}^{\a} :=\hat D_\b{}^\a-D_\b{}^\a,
\quad F^\a :=\hat E^\a-E^\a,
\quad A :=\hat B-B
\end{equation*}
and $\hat D_\b{}^\a, \hat E^\a, \hat B$ are the analogues for $\hat M$
of the functions \eqref{eq-CRherm2}
restricted to $M$.

Let us first compute $C_\b{}^\a$ using the
pseudohermitian Gauss equation obtained in \S\ref{sec-gauss}.
In view of \eqref{eq-CRherm2}, we have
\begin{equation}\Label{Cab-start}
C_{\a\bar\b}=  \frac{i\hat
R_{\a\bar\b}}{{\hat n}+2}-\frac{i\hat Rg_{\a\bar\b}}
{2({\hat n}+1)({\hat n}+2)}-\frac{iR_{\a\bar\b}}{n+2}+\frac{iRg_{\a\bar\b}}{2(n+1)(n+2)},
\end{equation}
where the pseudohermitian Ricci curvature tensors
$\hat R_{\a\bar\b}$ and $R_{\a\bar\b}$ are related by
\begin{equation}\Label{eq-ricci}
R_{\a\bar\b} = \hat R_\mu{}^\mu{}_{\a\bar\b} - \omega_\mu{}^a{}_\a\,\omega^{\mu}{}_{a \bar \b}
= \hat R_{\a\bar\b} - \hat R_{a}{}^{a}{}_{\a\bar\b} - \omega_\mu{}^a{}_\a\,\omega^{\mu}{}_{a \bar \b},
\end{equation}
where we have used the contraction of \eqref{eq-relcurv}.
The pseudohermitian scalar curvature of $R$ is now obtained by
further contracting \eqref{eq-ricci}:
\begin{equation}\Label{eq-scalar}
R = \hat R_\mu{}^\mu - \hat R_{a}{}^{a}{}_{\mu}{}^\mu
- \o_\mu{}^a{}_\nu\, \o^{\mu}{}_{a}{}^\nu.
\end{equation}
Substituting (\ref{eq-ricci}--\ref{eq-scalar}) into \eqref{Cab-start},
we obtain
\begin{multline}\Label{cab}
C_{\a\bar\b}=
\frac{i\omega_\mu{}^a{}_\a\,\omega^\mu{}_a{}_{\bar\b}}{n+2}
+\frac{i\hat R_a{}^a{}_{\a\bar\b}}{n+2}
+ \Big(\frac{1}{{\hat n}+2}-\frac{1}{n+2}\Big)i\hat R_{\a\bar\b}\\
- \frac{i\omega_\mu{}^a{}_\nu\,\omega^\mu{}_a{}^{\nu}g_{\a\bar\b}}{2(n+1)(n+2)}
-\frac{i (\hat R_a{}^a{}_{\mu}{}^{\mu} +  \hat R_{\mu}{}^{\mu})\, g_{\a\bar\b}} {2(n+1)(n+2)}
-\frac{i\hat Rg_{\a\bar\b}}{2({\hat n}+1)({\hat n}+2)}.
\end{multline}
We next express the contractions of the (tangential)
pseudohermitian tensor of $\hat M$ appearing here
in terms of its (tangential) pseudoconformal curvature tensor
and the (pseudohermitian) Ricci and scalar curvature.
For this, we contract the analogue of (\ref{eq-curvid})
(using the analogue of \eqref{Ricci}) for $\hat M$
over $a=n+1,\ldots,\hat n$, and restrict as before to the situation where
$g_{A\bar B}=\delta_{A\bar B}$ (so that, in particular, $g_{a\bar \mu}=0$):
\begin{equation}\Label{eq-long3}
\hat R_{a}{}^a{}_{\a\bar\b}
=\hat S_{a}{}^a{}_{\a\bar\b}
+\frac{\hat R_a{}^a g_{\a\bar\b}
+({\hat n}-n)\hat R_{\a\bar\b}}
{{\hat n}+2}
-\frac{({\hat n}-n)\hat R\, g_{\a\bar\b}}
{({\hat n}+1)({\hat n}+2)}.
\end{equation}
Contracting again for $\mu=1,\ldots,n$, we further obtain
\begin{equation}\Label{eq-long4}
\hat R_{a}{}^a{}_{\mu}{}^\mu
=\hat S_{a}{}^a{}_{\mu}{}^\mu
+\frac{n\hat R_a{}^a +({\hat n}-n)\hat R_{\mu}{}^\mu}
{{\hat n}+2}
-\frac{n({\hat n}-n) \hat R}{({\hat n}+1)({\hat n}+2)}.
\end{equation}
After substituting (\ref{eq-long3}--\ref{eq-long4}) in \eqref{cab}
and simplifying, we obtain
\begin{equation}\Label{eq-Cab}
C_{\a\bar\b}
=\frac{i(\hat S_a{}^a{}_{\a\bar\b} + \omega_\mu{}^a{}_\a\,\omega^\mu{}_a{}_{\bar\b})}{n+2}
-\frac{i(\hat S_a{}^a{}_{\mu}{}^{\mu}
+ \omega_\mu{}^a{}_\nu\,\omega^\mu{}_a{}^{\nu})g_{\a\bar\b}}
{2(n+1)(n+2)}.
\end{equation}
Thus the coefficients $C_{\a\bar\b}$ are completely
determined by the second fundamental form $\omega_\a{}^a{}_\b$
and the tangential pseudoconformal curvature tensor
$\hat S_{A\bar B C \bar D}$ of $\hat M$ pulled back to $M$.

Let us proceed to determine the $F^\a$. We compute
$d\hat\phi^\a$ modulo $\theta$ by differentiating the second
equation in (\ref{eq-CRMhatM}) and using
the sixth structure equation in \eqref{eq-cmw} pulled back to $M$ and \eqref{pbvanish}:
\begin{equation}\Label{eq-dphia1}
d\hat\phi^\a
=\phi^\mu\wedge\phi_\mu{}^\a
-\frac{1}{2}\psi\wedge\theta^\a
+\Phi^\a
+(dC_\nu{}^\a-C_\mu{}^\a\o_\nu{}^\mu)\wedge\theta^\nu
+iF^\a g_{\mu\bar\nu}\theta^\mu\wedge\theta^{\bar\nu}\ \mod\theta.
\end{equation}
On the other hand, by substituting (\ref{eq-CRMhatM}) into the
structure equation for $d\hat\phi^\a$ and using Proposition~\ref{prop1}
and  \eqref{eq-secform} to replace $\hat \phi^a\wedge\hat\phi_a{}^\a$ by
$\hat D_{\mu}{}^a\theta^\mu\wedge\omega_a{}^\a \mod\theta$, we obtain
\begin{multline}\Label{eq-dphia2}
d\hat\phi^\a
=\phi^\mu\wedge\phi_\mu{}^\a
+C_\nu{}^\mu\theta^\nu\wedge\o_\mu{}^\a
+\hat D_{\mu}{}^a\theta^\mu\wedge\omega_a{}^\a\\
-\frac{1}{2}(\psi\wedge\theta^\a
+iF_\mu\theta^\mu\wedge\theta^\a
-iF_{\bar\nu}\theta^{\bar\nu}\wedge\theta^\a)
+\hat \Phi^\a
\mod\theta
\end{multline}
and hence, using the covariant derivative (see \S\ref{sec-prelim}),
\begin{multline}\Label{eq-Fa2}
\Phi^\a + \nabla C_{\nu}{}^\a\wedge\theta^\nu+iF^\a
g_{\mu\bar\nu}\theta^\mu\wedge \theta^{\bar\nu}
=\hat D_\nu{}^a\theta^\nu\wedge\omega_a{}^\a\\
-\frac{1}{2}(iF_\mu\theta^\mu\wedge\theta^\a
+iF_{\bar\nu}\theta^\a\wedge\theta^{\bar\nu})
+\hat\Phi^\a\
\mod\theta.
\end{multline}
Identifying coefficients in front of
$\theta^\mu\wedge\theta^{\bar\nu}$
and using \eqref{eq-CRtor} we obtain:
\begin{equation}\Label{eq-FandD1}
\hat D_\mu{}^a\omega_a{}^\a{}_{\bar\nu} - iF^\a
g_{\mu\bar\nu}-\frac{1}{2}iF_{\bar\nu}\delta_\mu{}^\a=
V^\a{}_{\mu\bar\nu}-\hat V^\a{}_{\mu\bar\nu}-C_\mu{}^\a{}_{;\bar\nu}
\end{equation}
or, equivalently,
\begin{equation}\Label{eq-FandD2}
\hat D_\mu{}^a\omega_a{}^\a{}^\gamma - iF^\a
\delta_{\mu}{}^\gamma-\frac{1}{2}iF^\gamma\delta_\mu{}^\a=
V^\a{}_{\mu}{}^\gamma-\hat
V^\a{}_{\mu}{}^\gamma-C_\mu{}^\a{}_{;}{}^\gamma.
\end{equation}
Recall that we have the symmetry relation
$\omega_\a{}^a{}_\gamma=\omega_\gamma{}^a{}_\a$ by \eqref{eq-secform}
which implies
$\omega_a{}^{\a\gamma}=\omega_a{}^{\gamma\a}$ in view of \eqref{eq-consym2}.
We also have the symmetries
$V^\a{}_{\mu}{}^\gamma=V^\gamma{}_{\mu}{}^\a$ and $\hat
V^\gamma{}_{\mu}{}^\a=\hat V^\a{}_{\mu}{}^\gamma$ (see \cite[(A.14)]{CM}). By
subtracting (\ref{eq-FandD2}) from the same equation with $\gamma$
and $\a$ interchanged, we then obtain
\begin{equation}\Label{eq-Fonly}
\frac{1}{2}iF^\gamma\delta_\mu{}^\a-\frac{1}{2}iF^\a
\delta_{\mu}{}^\gamma
= C_\mu{}^\a{}_{;}{}^\gamma - C_\mu{}^\gamma{}_{;}{}^\a,
\end{equation}
which, after setting $\gamma=\mu$ and summing implies
\begin{equation}\Label{eq-Fonly1}
\frac{n-1}{2}\,iF^\a =
C_\mu{}^\mu{}_{;}{}^\a -C_\mu{}^\a{}_{;}{}^\mu.
\end{equation}
Thus, for $n\geq 2$, the coefficients $F^\a$
(with respect to a fixed admissible coframe $(\theta,\theta^\a)$)
are completely determined by the covariant derivatives of $C_\a{}^\b$
(given in terms of the uniquely defined forms $\o_\a{}^\b$)
and hence, in view of \eqref{eq-Cab},
by $\omega_\a{}^a{}_\b$ and $\hat S_{A\bar B C \bar D}$.
It now follows from the first two equations in \eqref{eq-CRMhatM}
that the pulled back forms $\hat\phi_\a{}^\b$ and $\hat\phi^\a$ to $M$
are completely determined
in terms of the fixed data of $M$ and the pulled back
tangential pseudoconformal curvature tensor $\hat M$.

Next, we repeat the procedure above to $\hat\psi$.
By differentiating the third equation in \eqref{eq-CRMhatM}
and using the last equation of \eqref{eq-cmw}
and its analogue for $\hat M$ pulled back to $M$, we obtain
\begin{equation}\Label{}
2i\hat\phi^\mu\wedge\hat\phi_\mu + 2i\hat\phi^a\wedge\hat\phi_a + \hat\Psi
= 2i\phi^\mu\wedge\phi_\mu + \Psi + d(iF_\mu\theta^\mu-iF_{\bar\nu}\theta^{\bar\nu} +A\theta)
\end{equation}
and using again (\ref{eq-CRMhatM}),
\begin{multline}\Label{eq-AandD}
2i(-C_{\mu}{}^\gamma \phi_\gamma\wedge \theta^\mu
+C_{\bar\nu\gamma}\phi^\gamma\wedge\theta^{\bar\nu})
+2i ( C_\mu{}^\gamma C_{\bar\nu\gamma}+ \hat D_\mu{}^a\hat D_{\bar\nu a})
\theta^\mu\wedge\theta^{\bar\nu}
+ \hat\Psi\\
= \Psi+i\nabla F_\mu\wedge\theta^\mu
-iDF_{\bar\nu}\wedge\theta^{\bar\nu}
+i Ag_{\mu\bar\nu}\theta^\mu\wedge\theta^{\bar\nu} \mod \theta.
\end{multline}
By identifying the coefficients in front of
$\theta^\mu\wedge\theta^{\bar\nu}$, using (\ref{eq-CRtor})
and (\ref{eq-CRherm1}), we obtain
\begin{multline}\Label{eq-AandD2}
2i(C_\mu{}^\gamma D_{\bar\nu\gamma}+C_{\bar\nu\gamma}D_\mu{}^\gamma)
+2i(C_\mu{}^\gamma C_{\bar\nu\gamma}+\hat D_\mu{}^a\hat D_{\bar\nu a})
-2i\hat P_{\mu\bar\nu}\\
=-2iP_{\mu\bar\nu}-iF_{\mu;\bar\nu}-iF_{\bar\nu;\mu}+iAg_{\mu\bar\nu}.
\end{multline}
By contracting, we obtain
\begin{equation}\Label{eq-AandD3}
nA =(F_{\mu;}{}^\mu+F^\mu{}_{;\mu})-2 \hat P_{\mu}{}^\mu
+2 (C_\mu{}^\gamma C^\mu{}_{\gamma} + C_\mu{}^\gamma D^\mu{}_{\gamma} + C^\mu{}_\gamma D_\mu{}^\gamma
+\hat D_\mu{}^a\hat D^\mu{}_{a}),
\end{equation}
where we have used the trace condition $P_\mu{}^\mu=0$.

The coefficients $D^\mu{}_{\gamma}$ are given by \eqref{eq-CRherm2}
and hence are determined by the coframe $(\theta,\theta^\a)$.
Suppose now in addition that the second fundamental form $(\o_\a{}^a{}_\b)$
is {\em nondegenerate} in the sense that the collection of vectors
$b_{\a\b}:=(\omega_\a{}^{n+1}{}_{\b},\ldots,\omega_\a{}^{\hat n}{}_{\b})$
spans $\bC^{{\hat n}-n}$.
Then $\o_\a{}^a{}_\b=\o_{\a \bar a \b}=-\o_{\bar a \a\b} = -\1{\o_{a}{}^{\a}{}_{\bar\b}}$
in view of \eqref{eq-consym2} and hence we can solve for the coefficients
$\hat D_\mu{}^a$ in (\ref{eq-FandD1}). Using
equation (\ref{eq-AandD3}), we can then determine $A$. Hence, all the
pulled back forms in (\ref{eq-CRMhatM}), as well as the forms
$\hat\phi_{\mu\bar a} = -\hat\phi_{\bar a\mu}$ (see the analogue of \eqref{eq-CRherm1} for $\hat M$),
are uniquely determined by the second
fundamental form and the pulled back pseudoconformal curvature of $\hat
M$. Now, if the manifold $\hat M$ is the sphere, its pseudoconformal curvature
vanishes identically (see \cite{CM}) and the forms above are {\it uniquely
determined by the second fundamental form}. We shall show that,
when $\hat M$ is the sphere, all the pulled back forms
$\hat\phi_A{}^B,\hat \phi^A,\hat\psi$ on $M$ are determined by the
second fundamental form. It remains to determine $\hat\phi_b{}^a$
and $\hat\phi^a$ in this case.
The choice of the 1-forms $\hat \phi_b{}^a$ on $M$ is always determined up to
\begin{equation}\Label{eq-indet}
\tilde{\hat\phi}_b{}^a=
\hat\phi_b{}^a+n_b{}^a{}_\gamma\theta^\gamma
+m_b{}^a{}_{\bar\gamma}\theta^{\bar\gamma}
+p_b{}^a\theta,\end{equation}
where, since $\hat \phi_{b\bar a}=-\hat \phi_{\bar a b}$, we must have
the symmetries \begin{equation}\Label{eq-symme} n_{b\bar
a\gamma}+\overline{m_{a\bar b\bar\gamma}}=p_{b\bar
a}+\overline{p_{a\bar b}}=0
\end{equation}
in view of the first equation in (\ref{eq-cmw}).
By the analogue of Proposition \ref{prop1} for $\hat M$
and (\ref{eq-secform}), we have
\begin{equation}\Label{eq-hatphia}
\hat\phi^a=\hat D_\b{}^a\theta^\b+\hat
E^a\theta,
\end{equation}
where we have used the vanishing of $\hat\theta^a$ on $M$.
Since $\hat D_\b{}^a$ has
already been determined, $\hat\phi^a$ is determined up to
\begin{equation}\Label{eq-hatphia2}
\tilde{\hat\phi}^a=\hat \phi^a + G^a\theta
\end{equation}
for some $G^a$.
Since every possible choice of $\hat\phi_b{}^a$ and $\hat\phi^a$ must satisfy the
structure equation \eqref{eq-cmw} for $d\hat\phi_\b{}^a$, we conclude, by
subtracting two such possibilities from each other and using
the analogue of \eqref{eq-CRherm1} for $\hat M$, \eqref{eq-secform} and the
vanishing of $\hat \Phi_\b{}^a$ for the sphere, that
\begin{equation}\Label{eq-hatphiab}
0=\big(\omega_\b{}^b{}_\a\theta^\a +  \hat D_\b{}^b\theta \big)\wedge
\big(\tilde {\hat\phi}_b{}^a-\hat\phi_b{}^a\big)+
ig_{\b\bar\mu} G^a \theta^{\bar\mu}\wedge\theta.
 \end{equation}
 By substituting (\ref{eq-indet}) in (\ref{eq-hatphiab}), we first
 conclude, using the nondegeneracy of $\omega_\a{}^a{}_\b$,  that
 $m_b{}^a{}_{\bar\b}=0$ and, hence, by (\ref{eq-symme}) that
 $n_b{}^a{}_\gamma=0$. Now, it easily follows from whatever is
 left in (\ref{eq-hatphiab}) that $p_b{}^a=G^a=0$. This completes
 the proof of the following result.

\begin{Thm}\Label{thm-conn}
Let $f\colon M\to \hat M$ be a CR-embedding of
a strictly pseudoconvex CR-manifold $M$ of dimension $2n+1$, $n\geq 2$,
into the unit sphere $\hat M=\5S$ in $\bC^{{\hat n}+1}$.
Denote by $(\omega_\a{}^a{}_\b)$ the second fundamental form of
$f$ relative to some coframe $(\theta,\theta^A)$ on $\hat M$
adapted to $f(M)$. If $(\omega_\a{}^a{}_\b)$ is nondegenerate,
then it uniquely determines $(\hat\phi_B{}^A,\hat \phi^A,\hat \psi)$.
That is, if there is another CR-embeddings $\tilde f\colon M\to \hat M$
having the same second fundamental form $(\omega_\a{}^a{}_\b)$
relative to some coframe $(\tilde\theta,\tilde \theta^A)$ on $\hat M$
adapted to $\tilde f(M)$ with
$\tilde f^*(\tilde\theta,\tilde\theta^\a)=f^*(\theta,\theta^\a)$,
then necessarily
\begin{equation}\Label{eq-thmconn1}
\tilde f^*(\tilde{\hat\phi}_B{}^A, \tilde{\hat\phi}^A, \tilde{\hat\psi}) =
f^*(\hat\phi_B{}^A, \hat\phi^A, \hat\psi),
\end{equation}
where $(\hat\phi_B{}^A,\hat\phi^A,\hat\psi)$
and $(\tilde{\hat\phi}_B{}^A,\tilde{\hat\phi}^A,\tilde{\hat\psi})$
are the corresponding parts of parallellisms pulled back to $\hat M$.
\end{Thm}

\section{Finitely nondegenerate mappings}\Label{sec-degsff}

In Theorem \ref{thm-conn} above, we assumed that the second
fundamental form of the embedding is nondegenerate or,
equivalently, that the embedding is $2$-nondegenerate. This
condition can be relaxed to merely requiring that the embedding is
{\it finitely nondegenerate}, i.e. $k$-nondegenerate for some $k$.
To do this, we first interpret $k$-nondegeneracy of the
embedding as a condition on the second fundamental form
that we view as a section in the vector bundle of $\C$-bilinear maps
$$T^{1,0}_pM\times T^{1,0}_pM\to T^{1,0}_p\hat M/T^{1,0}_pM, \quad p\in M.$$
For sections of this bundle we have the covariant differential
induced by the pseudohermitian connections $\nabla$ and $\hat\nabla$
on $M$ and $\hat M$ respectively:
\begin{equation}\Label{o-der}
\nabla\omega_{\a}{}^a{}_\b
=d\omega_\a{}^a{}_\b
-\omega_\mu{}^a{}_\b\,\hat\o_\a{}^\mu
+\omega_\a{}^b{}_\b\,\hat\o_b{}^a-
\omega_\a{}^a{}_\mu\,\hat\o_\b{}^\mu,
\end{equation}
As above, we also use e.g.\
$\o_\a{}^a{}_{\b;\gamma}$ to denote its component
in the direction $\theta^\gamma$.
Higher order covariant derivatives
$\o_\a{}^a{}_{\b;\gamma_1,\ldots,\gamma_l}$
are defined inductively in a similar way:
\begin{equation}\Label{eq-cov1}
\nabla\omega_{\gamma_1}{}^a{}_{\gamma_2;\gamma_3\ldots\gamma_j}=
d\omega_{\gamma_1}{}^a{}_{\gamma_2;\gamma_3\ldots\gamma_j}+
\omega_{\gamma_1}{}^b{}_{\gamma_2;\gamma_3\ldots\gamma_j} \, \hat\o_b{}^a
-\sum_{l=1}^j
\omega_{\gamma_1}{}^a{}_{\gamma_2;\gamma_3\ldots\gamma_{l-1}\,\mu\,\gamma_{l+1}\ldots\gamma_j}
\,\hat\o_{\gamma_l}{}^\mu,
\end{equation}

To obtain the desired relation between finite nondegeneracy
and covariant derivatives of $\o_{\a}{}^a{}_\b$,
we begin by applying $L_\gamma$ to \eqref{eq-omega}
and using the analogue of \eqref{eq-consymmetry} for $\hat M$:
\begin{equation}\Label{eq-eq2}
L_\gamma L_\a L_\b (\hat\rho_{\bar Z'}\circ f)
= L_{\gamma}\lrcorner d(\o_{\a\bar a\b} \theta^{\bar a})
= L_{\gamma}\lrcorner (d\o_{\a\bar a\b}\wedge \theta^{\bar a}
- \o_{\a\bar a\b} \,\hat\o_{\bar b}{}^{\bar a} \wedge  \theta^{\bar b} )
\mod \theta,\theta^{\bar\a}.
\end{equation}
We now observe that the expression in the brackets on the right hand side
coincides with $\nabla\o_{\a\bar a\b}\wedge \theta^{\bar a}$
modulo $\bar E_2$ (as defined in \S\ref{sec-main}) which is spanned
by the conjugated forms \eqref{e2-forms}.
Hence we obtain
\begin{equation*}
L_{\gamma} L_\a L_\b (\hat\rho_{\bar Z'}\circ f) =\o_{\a\bar
a\b;\gamma}\,\theta^{\bar a} \ \mod \bar E_2
\end{equation*}
and, by induction,
\begin{equation*}
L_{\gamma_l}\ldots L_{\gamma_1} L_\a L_\b (\hat\rho_{\bar Z'}\circ
f) = \o_{\a\bar a\b;\gamma_1,\ldots,\gamma_l}\,\theta^{\bar a} \
\mod \bar E_{l+1}
\end{equation*}
for any integer $l > 1$, where $E_{l+1}$ is as defined in \S\ref{sec-main}
and is spanned by
\begin{equation}\Label{el-span}
\theta,\theta^{\bar a},\o_{\a\bar a\b;\gamma_1,\ldots,\gamma_{s}}
\,\theta^{\bar a}, \quad 0\le s\le l-1.
\end{equation}
Hence, we obtain the following generalization
of the equivalence of $2$-nondegeneracy with nondegeneracy
of the second fundamental form:

\begin{Pro}\Label{pro-basic}
For $k_0\ge 2$, the embedding $f\colon M\to\hat M$ is $k_0$-nondegenerate at
a point $p$ if and only if the equality
\begin{equation}\Label{span-nondeg}
\span \{\o_{\gamma_1}{}^a{}_{\gamma_2;\gamma_3,\ldots,\gamma_{l}} \, L_a,
\; 2\le l\le k \} (p) = N_{f(p)}M
\end{equation}
holds for $k\ge k_0$ but not for $k<k_0$.
Furthermore, $f$ is $(k_0,s_0)$-nondegenerate at $p$ if and only if
the codimension of the span on the left hand side in \eqref{span-nondeg}
in the normal space $N_{f(p)}M$ is equal to $s_0$ for $k\ge k_0$
and strictly smaller for $k<k_0$.
\end{Pro}

We shall now prove the following alternative to Theorem
\ref{thm-conn}. (The reader should note that Theorem
\ref{thm-conn2} is more general than Theorem \ref{thm-conn} in
that it allows for a higher degree $k$ of nondegeneracy, but is a
little weaker than Theorem \ref{thm-conn} in the special case
$k=2$, i.e. in the case of nondegenerate second fundamental form,
since the theorem below requires also the covariant derivatives
$\omega_{\a}{}^a{}_{\b;\gamma}$ to determine the induced
pseudoconformal connection.)

\begin{Thm}\Label{thm-conn2}
Let $f\colon M\to \hat M$ be a $k$-nondegenerate CR-embedding of
a strictly pseudoconvex CR-manifold $M$ of dimension $2n+1$, $n\geq 2$,
into the unit sphere $\hat M=\5S$ in $\bC^{{\hat n}+1}$.
Denote by $(\omega_\a{}^a{}_\b)$ the second fundamental form of
$f$ relative to some coframe $(\theta,\theta^A)$ on $\hat M$
adapted to $f(M)$. Then the covariant derivatives
$(\o_\a{}^a{}_{\b;\gamma_1\ldots,\gamma_j})_{0\leq j\leq k-1}$ uniquely
determine $(\hat\phi_B{}^A,\hat\phi^A,\hat\psi)$.
That is, if $\tilde f\colon M\to\hat M$ is another CR-embedding
whose second fundamental form with respect to
a coframe $(\tilde\theta,\tilde \theta^A)$ on $\hat M$
adapted to $\tilde f(M)$ with $\tilde f^*(\tilde\theta,\tilde\theta^\a)=f^*(\theta,\theta^\a)$
is denoted by ${\tilde\omega}_\a{}^a{}_\b$
then
\begin{equation}\Label{eq-thmconn20}
\omega_\a{}^a{}_{\b;\gamma_1\ldots\gamma_j}
={\tilde\omega}_\a{}^a{}_{\b;\gamma_1\ldots\gamma_j},\quad \, 0\leq j\leq k-1,
\end{equation}
implies
$\tilde f^*(\tilde{\hat\phi}_B{}^A,\tilde{\hat\phi}^A,\tilde{\hat\psi}) = f^*(\hat\phi_B{}^A,\hat\phi^A,\hat\psi)$.
\end{Thm}

\begin{proof}
Let us for brevity use the notation $(\hat\phi_B{}^A,\hat \phi^A,\hat \psi)$
and $(\tilde{\hat\phi}_B{}^A,\tilde{\hat\phi}^A,\tilde{\hat\psi})$
also for the pullbacks to $M$.
If we examine the arguments above
leading  to Theorem \ref{thm-conn}, we see that the nondegeneracy
of the second fundamental form has not been used
to obtain $(\tilde{\hat\phi}_\b{}^\a,\tilde{\hat\phi}^\a)=(\hat\phi_\b{}^\a,\hat\phi^\a)$.
To compare the normal components $\hat\phi_b{}^a$ and $\tilde{\hat\phi}_b{}^a$
we write as before the general relation between them in the form (\ref{eq-indet}--\ref{eq-symme}).
In view of \eqref{eq-cov1}, the analogues of \eqref{eq-CRherm1} for $\hat M$
and the fact that $\hat\o_\b{}^\a=\o_\b{}^\a$ (see \S\ref{sec-subs}),
the equations \eqref{eq-thmconn20} imply that
\begin{equation}\Label{eq-normconn1}
\omega_{\gamma_1}{}^b{}_{\gamma_2;\gamma_3\ldots\gamma_j}n_b{}^a{}_\gamma=0,
\quad 2\leq j\leq k.
\end{equation}
Since we assume $f$ to be $k$-nondegenerate,
we conclude from Proposition~\ref{pro-basic} that
$n_b{}^a{}_\gamma=0$ and hence, by \eqref{eq-symme}, also $m_b{}^a{}_{\bar\gamma}=0$.
It follows that
\begin{equation}\Label{eq-normconn2}
\tilde{\hat\phi}_b{}^a=\hat\phi_b{}^a+p_b{}^a\theta
\end{equation}
which, in turn, implies $\tilde{\hat\o}_b{}^a=\hat\o_b{}^a \mod \theta$
and therefore
\begin{equation}\Label{eq-conjconn}
\o_\a{}^a{}_{\b;\bar\gamma}=\tilde\o_\a{}^a{}_{\b;\bar\gamma}.
\end{equation}
By the analogue of \eqref{eq-CRherm1} for $\hat M$ and \eqref{eq-secform},
\begin{equation}\Label{eq-phialphaa}
\omega_\a{}^a{}_\b\theta^\b+\hat D_\a{}^a\theta = \hat\phi_\a{}^a,
\end{equation}
and similarly for $\tilde{\hat\phi}_\a{}^a$. We claim that
$\hat D_\a{}^a=\tilde{\hat D}_\a{}^a$. To see this, we
differentiate \eqref{eq-phialphaa},
use the fourth identity in \eqref{eq-cmw} for $d\hat\phi_\a{}^a$
and identify the coefficients of $\theta^\b\wedge\theta^{\bar\gamma}$:
\begin{equation}\Label{eq-codazzi2}
\omega_\a{}^a{}_{\b;\bar\gamma}=i(g_{\a\bar\gamma}\hat
D_\b{}^a + g_{\b\bar\gamma}\hat D_\a{}^a),
\end{equation}
whose contraction yields
\begin{equation}\Label{eq-cod3}
\omega_\mu{}^a{}_{\b;}{}^\mu=i(n+1)\hat D_\b{}^a.
\end{equation}
Thus, $\hat D_\a{}^a$ is uniquely determined by
$\omega_\a{}^a{}_{\b,\bar\gamma}$ proving the claim in view
of \eqref{eq-conjconn}.
Then $\tilde{\hat\phi}_\a{}^a={\hat\phi}_\a{}^a$ by \eqref{eq-phialphaa}
and the analogue of \eqref{eq-CRherm1} for $\hat M$ implies
in view of \eqref{eq-secform} that  $\tilde{\hat\phi}^a$ and $\hat\phi^a$
are related by \eqref{eq-hatphia2} with suitable $G^a$.
Moreover, it follows from equation \eqref{eq-AandD3} that
$A=\tilde A$ and hence $\hat\psi=\tilde{\hat\psi}$.

To finish the proof, we must show that $p_a{}^b= G^a=0$.
For this, we rewrite \eqref{eq-hatphiab}
\begin{equation}\Label{eq-end1}
0=\o_\a{}^b{}_\b \, p_b{}^a\theta^\a\wedge\theta
+ig_{\b\bar\nu} G^a\theta^{\bar\nu}\wedge\theta
\end{equation}
from which it follows that $G^a=0$ (and $\omega_\a{}^a{}_\b \, p_a{}^b=0$,
which does not necessarily imply $p_a{}^b=0$ unless
the second fundamental form is nondegenerate). Hence $\hat\phi^a=\tilde{\hat\phi}^a$.
To conclude that $p_a{}^b=0$, we
differentiate the equation \eqref{eq-normconn2}, use the structure equations
for $\hat\phi_a{}^b$ and $\tilde{\hat\phi}_a{}^b$ and identify the coefficient in front
of $\theta^\a\wedge\theta^{\bar\b}$ that yields $ig_{\a\bar\b}p_a{}^b=0$
and hence clearly implies $p_a{}^b=0$.
The proof of Theorem \ref{thm-conn2} is complete.
\end{proof}

Theorem \ref{thm-conn2} implies that the induced pseudoconformal
connection is determined by the second fundamental form and its
covariant derivatives with respect to this connection. We shall
now prove these covariant derivatives, in turn, are controlled by
the (pseudoconformal) Gauss equation. The result is the following.

\begin{Thm}\Label{thm-secform2}
Let $f\colon M\to \hat M$ be a CR-embedding of a strictly pseudoconvex CR
manifold of dimension $2n+1$ into the unit sphere $\hat M=\5S$
in $\bC^{{\hat n}+1}$ with ${\hat n}-n<n/2$.
Denote by $(\o_\a{}^a{}_\b)$
the second fundamental form of $f$ relative
to a coframe $(\theta,\theta^A)$ on $\hat M$ adapted to $f(M)$
and by $\o_\a{}^a{}_{\b;\gamma_1,\ldots,\gamma_l}$ its
covariant derivatives $($with respect to the pseudohermitian structures on $M$ and $\hat M)$.
Fix an integer $k\geq 2$ such that the space $E_k=E_k(p)$
$($spanned by \eqref{el-span}$)$ is of constant dimension for $p\in M$.
Then the collection of derivatives $\o_{\gamma_1}{}^a{}_{\gamma_2;\gamma_3,\ldots,\gamma_l}$
for $2\le l\le k$ is uniquely determined by $(\theta,\theta^\a)$
up to unitary transformations of
$\C^{\hat n-n}$,
i.e.\ for any other CR-embedding $\tilde f\colon M\to \hat M$ and any coframe
$(\tilde\theta,\tilde \theta^A)$ on $\hat M$ adapted to $\tilde f(M)$ with
$\tilde f^*(\tilde\theta,\tilde\theta^\a) = f^*(\theta,\theta^\a)$,
one has
\begin{equation}\Label{cov-eq}
\tilde{\o}_{\gamma_1}{}^a{}_{\gamma_2;\gamma_3,\ldots,\gamma_l}
= \o_{\gamma_1}{}^a{}_{\gamma_2;\gamma_3,\ldots,\gamma_l},
\quad 2\le l\le k,
\end{equation}
after  a possible unitary change of $(\tilde\theta^a)$ $($near any given point$)$.
\end{Thm}

\begin{proof}
As in the proof of Lemma~\ref{alg} and Corollary~\ref{thm-secform},
it is sufficient to prove that the collections of vectors
\begin{equation}\Label{collections}
\big(\o_{\gamma_1}{}^a{}_{\gamma_2;\gamma_3,\ldots,\gamma_l} L_a\big)_{2\le l\le k}
\quad\text{ and }\quad \big(\tilde\o_{\gamma_1}{}^a{}_{\gamma_2;\gamma_3,\ldots,\gamma_l}
\tilde L_a\big)_{2\le l\le k}
\end{equation}
have the same scalar products with respect to $(g_{a\bar b})$
(where $(T,L_A)$ and $(\tilde T,\tilde L_A)$ are dual frames to $(\theta,\theta^A)$
and $(\tilde\theta,\tilde\theta^A)$ respectively), i.e.\ that
\begin{equation}\Label{scalar-products}
g_{a\bar b}\, \o_{\gamma_1}{}^a{}_{\gamma_2;\gamma_3,\ldots,\gamma_l}
\, \o_{\bar\delta_1}{}^{\bar b}{}_{\bar\delta_2;\bar\delta_3,\ldots,\bar\delta_s}
= g_{a\bar b}\, \tilde\o_{\gamma_1}{}^a{}_{\gamma_2;\gamma_3,\ldots,\gamma_l}
\, \tilde\o_{\bar\delta_1}{}^{\bar b}{}_{\bar\delta_2;\bar\delta_3,\ldots,\bar\delta_s},
\quad 2\le l,s\le k,
\end{equation}
(the constant dimension assumption on $E_k$ is needed to guarantee that that the unitary
change of $(\theta^a)$ can be made smooth; otherwise a nonsmooth change is still possible).
For this, we observe that $\hat S_{\a\bar\b\mu\bar\nu} = 0$ since $\hat M$ is the unit sphere
and rewrite the pseudoconformal Gauss equation \eqref{gauss} in this case as
\begin{equation}\Label{full-gauss}
0= S_{\a\bar \b\mu\bar\nu} + g_{a\bar b}\,\omega_\a{}^a{}_\mu\, \omega_{\bar\b}{}^{\bar b}{}_{\bar\nu}
+ T_{\a\bar\b\mu\bar\nu},
\end{equation}
where $T_{\a\bar\b\mu\bar\nu}$ is a tensor conformally equivalent to $0$
(see \S\ref{sec-prelim}),
i.e.\ a linear combination of $g_{\a\bar\b}$, $g_{\a\bar\nu}$, $g_{\mu\bar\b}$ and $g_{\mu\bar\nu}$.
More generally, we call a tensor
$T_{\a_1,\ldots,\a_r,\bar\b_1,\ldots,\bar\b_s}{}^{a_1\ldots a_t\bar b_1\ldots\bar b_q}$
(and similar tensors with different order of indices)
{\em conformally equivalent to $0$} or {\em conformally flat} if it is a linear combination of $g_{\a_i\bar\b_j}$
for $i=1,\ldots,r$, $j=1,\ldots,s$.

We will show \eqref{scalar-products} inductively by comparing covariant derivatives of
the right hand side in \eqref{full-gauss} and in its analogue for $(\tilde\theta,\tilde\theta^\a)$.
When calculating the derivatives, the following two remarks with be of importance.
First, since we have assumed the Levi form matrix $(g_{\a\bar\b})$ to be constant,
it is easy to see that covariant derivatives of conformally flat tensors (with respect to any connection)
are always conformally flat. Second, in view of \eqref{eq-codazzi2},
the ``mixed'' derivatives $\o_{\a}{}^{a}{}_{\b;\bar\gamma}$ are conformally flat
together with all their higher order derivatives.
Finally, in order to treat e.g.\ derivatives $\o_{\a}{}^{a}{}_{\b;\gamma\bar\delta}$
we shall need the following lemma, which describes how covariant derivatives commute.

\begin{Lem}\Label{lem-covdercom}
In the setting of Theorem~{\rm\ref{thm-secform2}}, for any $s\geq 2$, we have a relation
\begin{equation}\Label{eq-covdercom}
\omega_{\gamma_1}{}^a{}_{\gamma_2;\gamma_3\ldots\gamma_s\a\bar\b}\, -\,
\omega_{\gamma_1}{}^a{}_{\gamma_2;\gamma_3\ldots\gamma_s\bar\b\a}
= C_{\gamma_1\ldots\gamma_s\a\bar\b}{}^{\mu_1\ldots\mu_s}
\, \omega_{\mu_1}{}^a{}_{\mu_2;\mu_3\ldots\mu_s}
+ \,T_{\gamma_1\ldots\gamma_s\a\bar\b}{}^a,
\end{equation}
where the tensor $C_{\gamma_1\ldots\gamma_s\a\bar\b}{}^{\mu_1\ldots\mu_s}$
depends only on $(\theta,\theta^\a)$ $($and not on the embedding $f${}$)$
and $T_{\gamma_1\ldots\gamma_s\a\bar\b}{}^a$ is conformally flat.
\end{Lem}

\begin{proof}
By observing that the left hand side of the identity
\eqref{eq-covdercom} is a tensor,
it is enough to show, for each fixed point $p\in M$,
the identity at $p$ with respect to any particular
choice of $(\theta,\theta^A)$ depending smoothly on $p$.
By making a suitable unitary change of coframe $\theta^\a\to
u_\b{}^\a \theta^\b$ and $\theta^a\to u_b{}^a \theta^b$
(in the tangential and normal directions respectively), we may
chose a coframe smoothly depending on $p$
with $\hat\o_\a{}^\b(p)=\hat\o_a{}^b(p)=0$.
Relative to such a coframe, the left hand side of \eqref{eq-covdercom}
evaluated at $p$ equals, in view of \eqref{eq-cov1} and \eqref{eq-CRherm1},
modulo a conformally flat tensor,
to the coefficient in front of $\theta^\a\wedge\theta^{\bar\b}$ in the expression
\begin{equation}\Label{eq-exp}
\sum_{j=1}^s \omega_{\gamma_1}{}^a{}_{\gamma_2;\gamma_3\ldots\gamma_{j-1}\,\mu\,\gamma_{j+1}\gamma_s}\,
d\o_{\gamma_j}{}^\mu
- \omega_{\gamma_1}{}^b{}_{\gamma_2;\gamma_3\ldots\gamma_s}\,
d\hat\phi_b{}^a,
\end{equation}
where we have used that $\hat\o_{\a}{}^{\b}=\o_{\a}{}^{\b}$ (see \S\ref{sec-subs}).
Next, we observe, by examining the equations \eqref{eq-CRMhatM}
and \eqref{eq-Cab}, that the coefficients in front of
$\theta^\a$ and $\theta^{\bar\b}$ in the pulled back forms
$\hat\phi^\gamma$ are uniquely determined by $(\theta,\theta^\a)$
and the scalar products $g_{a\bar b}\,\omega_\a{}^a{}_\mu\,\omega_{\bar\b}{}^{\bar b}{}_{\bar\nu}$.
The latter are, in turn, uniquely determined by $(\theta,\theta^\a)$
in view of Lemma~\ref{alg} (ii) (applied to \eqref{full-gauss}).
We then compute $d\hat\phi_b{}^a$ using the structure equation
and the vanishing of $\theta^a$ on $M$ and of $\hat\phi_\b{}^\a$ and $\hat\phi_b{}^a$ at $p$ modulo $\theta$:
\begin{equation}\Label{eq-dphi}
d\hat\phi_b{}^a = \hat\phi_b{}^\mu \wedge \hat\phi_\mu{}^a  - i\delta_b{}^a \hat\phi_\mu \wedge \theta^\mu  \mod \theta.
\end{equation}
In view of (\ref{eq-consym2} -- \ref{eq-secform}) and \eqref{eq-CRherm1},
the first term on the right hand side does not contribute to the the coefficient in front of
$\theta^\a\wedge\theta^{\bar\b}$. Hence we conclude that the coefficient in front of
$\theta^\a\wedge\theta^{\bar\b}$ in \eqref{eq-exp}, at $p$,
is of the required form.
This completes the proof of Lemma \ref{lem-covdercom}.
\end{proof}

We now return to the proof of Theorem \ref{thm-secform2}
and show the required indentities \eqref{scalar-products}
by induction on the number of indices $l+s$. The case $l+s=4$ (i.e.\ $l=s=2$)
follows from Lemma~\ref{alg} (ii) applied to \eqref{full-gauss}.
Since we have chosen coframes for which the Levi form
$g_{A\bar B}$ is constant, we have
$$
\hat\nabla g_{A\bar B}= dg_{A\bar B}-g_{C\bar B} \,\hat\o_{A}{}^C
-g_{A\bar D} \, \hat\o_{\bar B}{}^{\bar D} =
-(\hat\o_{A \bar B}+ \hat\o_{\bar B A}) = 0.
$$
In order to show \eqref{scalar-products} for any $2\le l,s\le k$ with $l+s>2$
we covariantly differentiate \eqref{full-gauss} for $\a=\gamma_1$, $\mu=\gamma_2$,
$\b=\delta_1$, $\nu=\delta_2$ with respect to the indices $\gamma_3,\ldots,\gamma_l$
and $\bar\delta_3,\ldots,\bar\delta_s$ consequently.
We obtain that the left hand side in \eqref{scalar-products}
is uniquely determined  by the covariant derivatives of
$S_{\a\bar\b\mu\bar\nu}$ modulo a conformally flat tensor
and modulo other scalar product terms involving
either derivatives of $\o_{\gamma_1}{}^a{}_{\gamma_2}$
with respect to indices one of which is $\delta_j$
or derivatives of $\o_{\bar\delta_1}{}^{\bar b}{}_{\bar\delta_2}$
with respect to indices one of which is $\gamma_j$.
If such an index appears at the beginning, the corresponding scalar product
is conformally flat in view of \eqref{eq-codazzi2}.
Otherwise, by interchanging indices and applying Lemma~\ref{lem-covdercom},
we see that the scalar product term is equal, again modulo a conformally flat tensor,
to a sum of scalar products of the expressions
on the right hand side of \eqref{eq-covdercom} and its derivatives.
These scalar products are linear combinations of
the expressions on the left hand side in \eqref{scalar-products}
with a smaller number of indices $l+k$.
Hence, by the induction, they are uniquely determined by $(\theta,\theta^\a)$.
Summarizing, we obtain the equality in \eqref{scalar-products}
modulo a conformally flat tensor
$T_{\gamma_1,\ldots,\gamma_l \bar\delta_1\ldots,\bar\delta_s}$.
The proof of Theorem~\ref{thm-secform2} is completed by applying \cite[Lemma 3.2]{Hu1} as in the proof of Lemma~\ref{alg}
to conclude that $T_{\gamma_1,\ldots,\gamma_l \bar\delta_1\ldots,\bar\delta_s}=0$.
\end{proof}

\section{Adapted $Q$-frames and proof of Theorem $\ref{thm-main3}$ in case $s=0$}\Label{sec-adapt}

We embed $\bC^{{\hat n}+1}$ in projective space $\bP^{{\hat n}+1}$ in the
usual way, i.e.\ as the set $\{\zeta^0\neq 0\}$ in the homogeneous
coordinates $[\zeta^0:\zeta^1:\ldots:\zeta^{{\hat n}+1}]$,
and, following \cite[\S1]{CM}, realize the unit sphere
as the quadric $Q$ given in  $\bP^{{\hat n}+1}$
by the equation $(\zeta,\zeta)=0$, where the Hermitian scalar product $(\cdot,\cdot)$
is defined by
\begin{equation}\Label{eq-scalarp}
(\zeta,\tau):=\delta_{A\bar B} \zeta^A\overline{\tau^B}
+\frac{i}{2}(\zeta^{\hat n+1}\1{\tau^0}-\zeta^0 \1{\tau^{\hat n+1}}).
\end{equation}
Recall here the index convention from previous sections:
capital Latin indices run from $1$ to ${\hat n}$.
A {\em $Q$-frame} (see e.g.\ \cite{CM})
is a unimodular basis $(Z_0,Z_1,\ldots,Z_{{\hat n}+1})$ of
$\bC^{{\hat n}+2}$ (i.e.\ $\det(Z_0,\ldots,Z_{{\hat n}+1})=1$)
such that $Z_0$ and $Z_{{\hat n}+1}$, as points in $\bP^{{\hat n}+1}$, are
on $Q$, the vectors $(Z_A)$ form an orthonormal basis
(relative to the inner product \eqref{eq-scalarp}) for the
complex tangent space to $S$ at $Z_0$ and $Z_{{\hat n}+1}$ and $(Z_{\hat n+1},Z_0)=i/2$.
(We shall denote the corresponding points in
$\bP^{{\hat n}+1}$ also by $Z_0$ and $Z_{{\hat n}+1}$; it should be clear from the
context whether the point is in $\bC^{{\hat n}+2}$ or $\bP^{{\hat n}+1}$.)
Equivalently, a $Q$-frame is any unimodular basis satisfying
$$(Z_A,Z_A) = 1, \quad (Z_{{\hat n}+1},Z_0)=-(Z_0,Z_{{\hat n}+1})=i/2,$$
while all other scalar products are zero.

On the space $\3B$ of all $Q$-frames there is a natural
free transitive action of the group $\SU({\hat n}+1,1)$ of unimodular
$({\hat n}+2)\times({\hat n}+2)$ matrices
that preserve the inner product \eqref{eq-scalarp}.
Hence, any fixed $S$-frame defines an isomorphism between $\3B$ and $\SU({\hat n}+1,1)$.
On the space $\3B$, there are
{\em Maurer-Cartan forms} $\pi_\Lambda{}^\Omega$, where capital Greek
indices run from $0$ to ${\hat n}+1$, defined by
\begin{equation}\Label{eq-MC}
dZ_\Lambda=\pi_\Lambda{}^\Omega Z_\Omega
\end{equation}
and satisfying $d\pi_\Lambda{}^\Omega = \pi_\Lambda{}^\Gamma\wedge\pi_\Gamma{}^\Omega$.
Here the natural $\C^{\hat n+2}$-valued $1$-forms $dZ_\Lambda$ on $\3B$ are defined
as differentials of the map $(Z_0,\ldots,Z_{\hat n+2})\mapsto Z_\Lambda$.

Recall \cite{CM,We79} that a smoothly varying $Q$-frame $(Z_\Lambda)=(Z_\Lambda(p))$ for $p\in Q$
is said to be {\em adapted} to $Q$ if $Z_0(p)=p$ as points of $\bP^{{\hat n}+1}$.
It is shown in \cite[\S5]{CM} that, if we use an adapted $Q$-frame to pull
back the 1-forms $\pi_\Lambda{}^\Omega$ from $\3B$ to $Q$ and set
\begin{equation}\Label{eq-speccof}
\theta:=\frac{1}{2}\pi_0{}^{{\hat n}+1}, \quad \theta^A:=\pi_0{}^A,
\quad \xi:=-\pi_0{}^0 + \1{\pi_0{}^0},
\end{equation}
we obtain a coframe $(\theta,\theta^A)$ on $Q$
and a form $\xi$ satisfying the structure equation
$$d\theta = i\delta_{A\bar B}\theta^A\wedge\theta^{\bar B} + \theta\wedge\xi.$$
In particular, it follows from \eqref{eq-speccof} that the coframe $(\theta^\a,2\theta)$
is dual to the frame defined by $(Z_A,Z_{\hat n+})$ on $M$
and hence depends only on the values of $(Z_\Lambda)$ at the same points.
Furthermore, there exists a unique section $M\to Y$ for which
the pullbacks of the forms $(\o,\o^\a,\phi)$ in \eqref{eq-CRpar}
are $(\theta,\theta^\a,\xi)$ respectively.
Then the pulled back forms $(\hat\phi_{B}{}^{A},\hat \phi^{A},\hat\psi)$
are given by \cite[(5.8b)]{CM}
\begin{equation}\Label{other-forms}
\hat\phi_{B}{}^{A}=\pi_B{}^A - \delta_B{}^A\pi_0{}^0,
\quad \hat\phi^A = 2\pi_{\hat n+1}{}^A,
\quad \hat\psi = -4\pi_{\hat n+1}^0.
\end{equation}
As in \cite[(5.30)]{CM}, the pulled back forms $\pi_\Lambda{}^\Omega$ can be
uniquely solved from (\ref{eq-speccof}--\ref{other-forms}):
\begin{equation}\Label{eq-piphi}
\begin{aligned}
(\hat n+2)\pi_0{}^0 &= - \hat\phi_C{}^C - \xi,&
\pi_0{}^A &= \theta^A,&
\pi_0{}^{{\hat n}+1} &= 2\theta,\\
\pi_A{}^0 &= -i\hat\phi_A,&
\pi_A{}^B &= \hat\phi_A{}^B+\delta_A{}^B\pi_0{}^0,&
\pi_A{}^{{\hat n}+1} &= 2i\theta_A,\\
4\pi_{{\hat n}+1}{}^0 &= -\hat\psi,&
2\pi_{{\hat n}+1}{}^A &= \hat\phi^A,&
(\hat n+2)\pi_{{\hat n}+1}{}^{{\hat n}+1} &= \hat\phi_{\bar D}{}^{\bar D} + \xi.
\end{aligned}
\end{equation}
Thus, the pullback of $\pi_\Lambda{}^\Omega$ is completely
determined by the pullbacks $(\theta,\theta^A,\xi,\hat\phi_{B}{}^{A},\hat \phi^{A},\hat\psi)$.

In an adapted $Q$-frame $(Z_\Lambda)$, the vector $Z_0$ giving the reference point is determined
up to a factor $t\in\C$ and a change
$(Z_0,Z_A,Z_{\hat n+1})\mapsto (tZ_0,Z_A,t^{-1}Z_{\hat n+1})$
results in the change $\theta\mapsto |t|^2\theta$ (see \cite[(5.10)]{CM}).
Following \cite{CM}, denote by $H_1$ the subgroup of $\SU({\hat n}+1,1)$
fixing the point $Z_0$ in $\bP^{\hat n+1}$ and the length of $Z_0$ in $\C^{\hat n+2}$.
Then, there is a surjective homomorphism $H_1\to G_1$,
explicitly described in \cite[(5.12)]{CM},
onto the group $G_1$ acting transitively on the coframes $(\o,\o^A,\o^{\bar A},\phi)$
(on the bundle $E$ where $\o$ is fixed) such that the relations (\ref{eq-piphi})
are preserved. We conclude, in particular, that
for any choice of an admissible coframe $(\theta,\theta^A)$ on $Q$
(as defined in \S\ref{sec-prelim}),
there exists an adapted $Q$-frame $(Z_\Lambda)$
such that (\ref{eq-piphi}) holds with $\xi=0$.

\begin{proof}[Proof of Theorem {\rm\ref{thm-main3}} in case $s=0$]
Let $f,\tilde f\colon M\to \5S$ be two CR-embeddings as in Theorem \ref{thm-main0}
and let $p\in M$ be a point where $M$ is $k_0$-nondegenerate.
Then the space $E_{k_0}(q)$ defined by \eqref{eq-ksdeg1}
is of maximal and hence constant dimension for $q\in M$ near $p$.
Choose any admissible coframe $(\theta,\theta^\a)$ on $M$ near $p$
(in the terminology of \S\ref{sec-prelim}).
Then Corollary~\ref{thm-admada} applied to $f$ and $\tilde f$
yields admissible coframes $(\hat\theta,\hat\theta^A)$ and
$(\tilde{\hat\theta},\tilde{\hat\theta}^A)$ on $\5S$
near $f(p)$ and $\tilde f(p)$ respectively
that are both adapted to $(\theta,\theta^\a)$
as defined after Corollary~\ref{thm-admada}.
Hence, by Theorem \ref{thm-secform2}, there exists a further unitary
change of the coframe $(\tilde{\hat\theta},\tilde{\hat\theta}^A)$ near $\tilde f(p)$
such that one has the equality \eqref{cov-eq} for the covariant derivatives
up to order $k_0-2$ of the second fundamental forms of $f$ and $\tilde f$.
Now Theorem~\ref{thm-conn2} can be applied to conclude
that the connection forms $(\hat\phi_B{}^A,\hat\phi^A,\hat\psi)$
and $(\tilde{\hat\phi}_B{}^A,\tilde{\hat\phi}^A,\tilde{\hat\psi})$
coincide when pulled back via $f$ and $\tilde f$ respectively.

As described above we realize the sphere $\5S$ as the quadric $Q$
and construct two adapted $Q$-frames $(Z_\Lambda)$ and $(\tilde Z_\Lambda)$
on $Q$ near $f(p)$ and $\tilde f(p)$ respectively
such that the relations \eqref{eq-piphi} and their analogues for $\tilde f$ hold.
We conclude that also the corresponding Maurer-Cartan forms $(\pi_\Lambda{}^\Omega)$
and $(\tilde\pi_\Lambda{}^\Omega)$ coincide when pulled back to $M$
via $f$ and $\tilde f$ respectively.

The end of the proof follows the line of \cite[Lemma (1.1)]{We78}.
Since the group $\SU(\hat n+1,1)$ acts transitively on
the space of all $Q$-frames, we can compose $\tilde f$
with such a transformation to obtain $Z_\Lambda(f(p))= \tilde Z_\Lambda(\tilde f(p))$.
Now the relation $f^*\pi_\Lambda{}^\Omega=\tilde f^*\tilde\pi_\Lambda{}^\Omega$
implies that, along any sufficiently small real curve in $M$, starting from $p$,
both $Q$-frames $Z_\Lambda(f(q))$ and $\tilde Z_\Lambda(\tilde f(q))$
satisfy the same first order ODE with the same initial values at $p$.
Hence we must have $Z_\Lambda(f(q))=\tilde Z_\Lambda(\tilde f(q))$ for $q\in M$ near $p$
and therefore $f(q)=\tilde f(q)$ since $Z_0(f(q))$ and $\tilde Z_0(f(q))$ give,
as points in $\bP^{\hat n+1}$, exactly the reference points
by the definition of the adapted frames.
Finally, the global coincidence $f\equiv \tilde f$ follows from the local one
by the well known facts that any CR-manifold of hypersurface type that is embeddable
into a sphere is automatically minimal and that that two CR-functions on a minimal connected CR-manifold
coincide if and only if they coincide in a nonempty open subset.
\end{proof}

\section{Degenerate CR-immersions and Proofs of Theorems~$\ref{thm-main2}$ and $\ref{thm-main3}$}\Label{sec-deg}

This section is devoted to the proof of Theorem~\ref{thm-main2},
where we assume that the CR-immersion $f\colon M\to \5S$ has degeneracy $s\ge 1$.

\begin{proof}[Proof of Theorem $\ref{thm-main2}$]
As in the proof of Theorem~\ref{thm-main3} in case $s=0$,
we choose an admissible coframe $(\theta,\theta^A)$ on $Q$ near $f(p_0)$, adapted to $f(M)$,
where $p_0\in M$ is now chosen such that the integer $s(p_0)$ defined in \eqref{deg1}
coincides with the degeneracy of $f$ (which is the minimum of $s(p)$ for all $p\in M$).
As before, denote by $(\omega_\a{}^a{}_\b)$
the second fundamental form of $f$ relative to this coframe.
Since the mapping $f$ is
constantly $(k_0,s_0)$-degenerate near $p_0$, by Proposition~\ref{pro-basic},
the dimension of the span in \eqref{span-nondeg} for $k=k_0$ is constant
and equals to $d:= {\hat n}-n-s_0$ for $p$ near $p_0$.
For the remainder of this proof, we shall use the indices $*,\#$ running over
the set $n+1,\ldots, n+d$ (possibly empty)
and the indices $i,j$ running over $n+d+1,\ldots, {\hat n}$.
Then after a unitary change of the $(\theta^a)$ if necessary, we may assume that
\begin{equation}\Label{eq-constdeg}
\span \{\o_{\gamma_1}{}^\#{}_{\gamma_2;\gamma_3,\ldots,\gamma_{l}} L_\#,
\; 2\le l\le k_0 \}= \span \{L_\#\},
\quad
\omega_{\gamma_1}{}^j{}_{\gamma_2;\gamma_3\ldots\gamma_l}\equiv 0,
\quad l\ge 2.
\end{equation}
We remark at this point that, if we can prove that the image $f(M)$
lies in a $({\hat n}+1-s)$-dimensional plane $P$, then clearly the
mapping $f\colon M\to \5S\cap P$ is $k_0$-nondegenerate since the
normal space of $M$ inside $\5S\cap P$ would be $\cong \bC^d$.
We can write
\begin{equation}\Label{eq-pound}
\hat\o_\#{}^j=\hat\o_\#{}^j{}_\mu\theta^\mu
+ \hat\o_\#{}^j{}_{\bar\nu}\theta^{\bar \nu}+\hat\o_\#{}^j{}_0\theta
\end{equation}
for suitable coefficients $\hat\o_\#{}^j{}_\mu$, $\hat\o_\#{}^j{}_{\bar\nu}$ and $\hat\o_\#{}^j{}_0$.
In view of the definition of the covariant derivatives (\ref{o-der}--\ref{eq-cov1}),
\eqref{eq-constdeg} immediately implies
\begin{equation}
\o_{\gamma_1}{}^\#{}_{\gamma_2;\gamma_3,\ldots,\gamma_{l}}\hat\o_\#{}^j{}_\mu=0
\end{equation}
for any $l\geq 2$, and hence,
\begin{equation}\Label{eq-H=0}
\hat\o_\#{}^j{}_\mu=0.
\end{equation}
Furthermore, in view of \eqref{eq-CRherm1} and \eqref{eq-secform}, \eqref{eq-constdeg} implies
\begin{equation}\Label{eq-phiia}
\hat\phi_\a{}^j=\hat D_\a{}^j\theta,
\quad \hat\phi^j=\hat D_\mu{}^j\theta^\mu+\hat E^j\theta
\end{equation} and
\begin{equation}\Label{eq-phiia2}
\hat\phi_\a{}^\#=\omega_\a{}^\#{}_\mu\theta^\mu
+\hat D_\a{}^\#\theta,
\quad \hat\phi^\#=\hat D_\mu{}^\#\theta^\mu
+\hat E^\#\theta.
\end{equation}
Differentiating $\hat\phi_\a{}^j$ and using the structure
equations and the vanishing of the pseudoconformal curvature of $\5S$,
we obtain
\begin{equation}\Label{eq-Dai}
\omega_\a{}^\#{}_\mu \theta^\mu\wedge
\hat\phi_\#{}^j
=i(g_{\a\bar\nu}\hat D_\mu{}^j
+ g_{\mu\bar\nu}\hat D_\a{}^j)
\theta^\mu\wedge\theta^{\bar\nu} \quad\mod\theta.
\end{equation}
By using  \eqref{eq-CRherm1} and \eqref{eq-pound} to compute $\hat\phi_\#{}^j$
and identifying the coefficients of $\theta^\mu \wedge \theta^{\bar\nu}$,
we conclude that
\begin{equation}\Label{eq-basphipi}
\omega_\a{}^\#{}_\mu \hat\o_\#{}^j{}_{\bar\nu}
=i(g_{\a\bar\nu}\hat D_\mu{}^j+
g_{\mu\bar\nu}\hat D_\a{}^j).
\end{equation}
Since the right hand side is conformally flat (see \S\ref{sec-degsff}),
$\#$ runs over $n+1,\ldots, n+d$, and $d= {\hat n}-n-s<n$ by the assumptions of
Theorem~\ref{thm-main2},
it follows from \cite[Lemma~3.2]{Hu1} (cf.\ e.g.\ the
proof of Lemma~\ref{alg}) that
\begin{equation}\Label{eq-omG}
\omega_\a{}^\#{}_\mu \hat\o_\#{}^j{}_{\bar\nu}=0
\end{equation} and, consequently, $\hat D_\a{}^j=0$.
Thus,
we have
\begin{equation}\Label{eq-have}
\hat\phi_\a{}^j=0, \quad \hat\phi^j=\hat E^j\theta.
\end{equation}
Substituting this information back into the structure equation for
$d\hat\phi_\a{}^j$, using $d\hat\phi_\a{}^j=0$,
\eqref{eq-pound},  \eqref{eq-phiia} and \eqref{eq-CRherm1},
and this time identifying terms involving $\theta^\mu\wedge\theta$, we
conclude that
\begin{equation}\Label{eq-omJ}
\omega_\a{}^\#{}_\mu (\hat\o_\#{}^j{}_0 + \hat D_\#{}^j)=0.
\end{equation}

We claim that $\hat\phi_\#{}^j=0$. To prove this it suffices to
show
\begin{equation}\Label{eq-assump}
\o_{\gamma_1}{}^\#{}_{\gamma_2;\gamma_3,\ldots,\gamma_{l}}\hat\o_\#{}^j{}_{\bar\nu}=
\o_{\gamma_1}{}^\#{}_{\gamma_2;\gamma_3,\ldots,\gamma_{l}}
(\hat\o_\#{}^j{}_0+\hat D_\#{}^j)=0,
\quad l\geq 2.
\end{equation}
We next differentiate \eqref{eq-pound} (with $\hat\o_\#{}^j{}_\mu=0$ by \eqref{eq-H=0}),
use \eqref{eq-CRherm1} and
the structure equation to compute $d\hat\phi_\#{}^j$ modulo $\theta$
and identify the coefficients in front of $\theta^\mu\wedge\theta^{\bar\nu}$
to conclude that the covariant derivative $\hat\o_\#{}^j{}_{\bar\nu;\gamma}$
is conformally flat. Here the covariant derivative is understood
with respect to the induced pseudohermitian connection
on the corresponding subspaces and quotient spaces, where the indices are running.
Hence, all higher order covariant derivatives $\hat\o_\#{}^j{}_{\bar\nu;\gamma_1,\ldots,\gamma_m}$
are also conformally flat (see \S\ref{sec-degsff}).
By taking covariant derivatives of \eqref{eq-omG}
we obtain that the first expression in \eqref{eq-assump}
is conformally flat. As before we conclude that it is actually zero
by \cite[Lemma~3.2]{Hu1} since $d<n$.
Similarly we obtain conformal flatness of
the covariant derivative $\hat\o_\#{}^j{}_{0;\gamma}$
by identifying the coefficients in front of
$\theta^\mu\wedge\theta$ in the same identity.
Again, all higher order derivatives are also conformally flat
and we obtain the vanishing of the second expression
in \eqref{eq-assump} by taking covariant derivatives
of \eqref{eq-omJ} and using \cite[Lemma~3.2]{Hu1}.
This proves \eqref{eq-assump} and hence the claim
$\hat\phi_\#{}^j=0$ in view of \eqref{eq-constdeg}.
Finally, substituting this information back into the structure
equation for $d\hat\phi_\a{}^j=0$ and using \eqref{eq-have},
 we conclude that $\hat\phi^j=0$.

Now, we are in a position to complete the proof of Theorem
\ref{thm-main2}. We have shown that
$\hat\phi_\a{}^j=\hat\phi_\#{}^j=\hat\phi^j=0$.
As in \S\ref{sec-adapt} we realize the sphere $\5S$ as the quadric $Q$
in $\bP^{\hat n+1}$ and choose an adapted $Q$-frame
$(Z_\Lambda)$ on $Q\cong\5S$ near $f(p)$.
We can further choose $(Z_\Lambda)$ satisfying \eqref{eq-piphi}
with $\xi=0$ for our given admissible coframe $(\theta,\theta^A)$
as described in \S\ref{sec-adapt}.
It follows then from the second row in (\ref{eq-piphi})
and the symmetry relation $\hat\phi_{A\bar B}=-\hat\phi_{\bar B A}$
that $\pi_{j}{}^\Omega=0$ for all $\Omega$ except possibly
$\Omega\in\{n+d+1,\ldots,\hat n\}$. Thus, we have
\begin{equation}\Label{eq-good}
dZ_i=\pi_i{}^j Z_j, \quad i,j\in \{n+d+1,\ldots,\hat n\},
\end{equation}
expressing the property that derivatives of the vectors $Z_i$
are linear combinations of $Z_j$ at every point.
It follows that the span of $Z_j=Z_j(p)$ must be constant in the Grassmanian
of $\C^{\hat n+2}$ for $p\in M$ near $p_0$.
Hence the vectors $Z_0(p)$ lie in the constant
$(\hat n+2-s)$-dimensional orthogonal subspace
(with respect to the inner product (\ref{eq-scalarp})).
Since, by definition of the adapted $Q$-frame,
$Z_0(p)$ gives the reference point $p$ in $\bP^{\hat n+1}$,
we conclude that $f(p)$ is contained in a fixed $(\hat n+1-s)$-dimensional
plane $P$ in $\C^{\hat n+1}$ for $p\in M$ near $p_0$.
The minimality of $M$ implies now the inclusion $f(p)\in P$ for all $p\in P$
by the uniqueness of CR-functions.
The last statement follows directly from the definitions.
The proof of Theorem \ref{thm-main2} is complete.
\end{proof}

\begin{proof}[Proof of Theorem $\ref{thm-main3}$ in general case]
Let  $f\colon M\to \5S$ be as in the theorem and denote by $s$
its degeneracy. By Theorem~\ref{thm-main2}, $f$ can be seen
as a $k_0$-nondegenerate CR-immersion near some point of $M$
into a smaller sphere $\5S'\subset\C^{n+k-s+1}$.
The degeneracy of the immersion obtained of $M$ into $\5S'$
is now $0$. The required conclusion follows now
from the case $s=0$ considered in \S\ref{sec-adapt}.
\end{proof}

\end{document}